\documentclass{amsart}
\usepackage{amssymb}
\usepackage{graphics}

\theoremstyle{plain}
\newtheorem{lemma}{Lemma}

\newtheorem{proposition}{Proposition}
\newtheorem{theorem}{Theorem}
\theoremstyle{definition}
\newtheorem{definition}{Definition}
\newtheorem{corollary}{Corollary}
\theoremstyle{definition}
\newtheorem{remark}{Remark}

\begin{document}
\title[The fundamental group of symplectic manifolds]{The fundamental group of symplectic manifolds with  Hamiltonian Lie group actions}
\author{Hui Li}$$
\address{Mathematics, University of Luxembourg\\
         162A, Ave de la Faiencerie\\
         L-1511, Luxembourg.}
\email{li.hui@uni.lu}
\begin{abstract}
   Let $(M, \omega)$ be a connected, compact symplectic manifold equipped with a Hamiltonian $G$ action, where $G$ is a  connected compact Lie group.
  Let $\phi$ be the moment map.
 In \cite{L}, we proved the following result for $G=S^1$ action:  as fundamental groups
   of topological spaces, $\pi_1(M)\cong\pi_1(M_{red})$,
   where $M_{red}$ is the symplectic quotient at any  value  of the moment map $\phi$, and
   $\cong$ denotes ``isomorphic to''. In this paper,
 we generalize this result to other connected compact Lie group $G$
 actions. We also prove that the above fundamental group is
 isomorphic to that of $M/G$. We briefly discuss the generalization
 of the first part of the results to non-compact manifolds with proper moment maps.
\end{abstract}
 \maketitle
 \section{Introduction}

  Let $(M, \omega)$ be a connected, compact symplectic manifold. Let us  assume a connected compact Lie group $G$ acts on $M$ in a Hamiltonian
 fashion with moment map  $\phi: M\rightarrow \mathfrak{g}^*$, where $\mathfrak{g}^*$ is the dual Lie algebra of $G$.
 Assume $\phi$ is equivariant with respect to the $G$ action, where $G$ acts on $\mathfrak{g}^*$ by the coadjoint action.
 Take a moment map value
 $a$ in $\mathfrak{g}^*$, the space $M_{G\cdot a}=\phi^{-1}(G\cdot a)/G$ is called the symplectic quotient
  or the reduced space at the coadjoint orbit $G\cdot a$. If $G_a$ is the stabilizer group of $a$ under
  the coadjoint
      action, by equivariance of the moment map,  the two reduced spaces are equal:
      $M_a=\phi^{-1}(a)/G_a=\phi^{-1}(G\cdot a)/G=M_{G\cdot a}$. We will use
      the two notations interchangeably.
 The space $M_{G\cdot a}$ can be a smooth symplectic manifold, or a symplectic orbifold, or a symplectic stratified space.
 The space $M_{G\cdot a}$, as a topological space,
 has a well defined fundamental group.
One has the notion of orbifold fundamental group which is different
 (see \cite{Ch} or \cite{TY} for
the definition of orbifold $\pi_1$, and see the example following Theorem~\ref{thm1}).
  In \cite{L}, we proved the following theorem:
 \begin{theorem}\label{thm1}
  Let $(M, \omega)$ be a connected, compact symplectic manifold equipped with a Hamiltonian $S^1$ action.
  Then, as fundamental groups of topological spaces, $\pi_1(M)\cong\pi_1(\hbox{minimum})\cong\pi_1(\hbox{maximum})\cong\pi_1(M_{red})$,
 where $M_{red}$ is the symplectic quotient at any  value in the image of the moment map $\phi$.
  \end{theorem}

  The proof mainly relies on Morse-Bott theory and symplectic reduction. It closely links the fundamental group
  of the ``manifold under a value $a$''
  $M^a=\{x\in M\mid \phi(x)\leq a\}$ with the fundamental group of the reduced spaces at values equal to or lower than
  $a$.

  The above theorem is not true for ``orbifold fundamental group''.
   For example,  let $S^1$ act on $(S^2\times S^2, 2\rho\oplus\rho)$ (where $\rho$ is the standard symplectic
   form on $S^2$) by $\lambda (z_1, z_2)=(\lambda^2 z_1, \lambda z_2)$. Let $0$  be the minimal value of the moment map. Then for $a\in (1, 2)$, $M_a$
   is an orbifold which is homeomorphic to $S^2$ and has two $\mathbb{Z}_2$ singularities. The orbifold $\pi_1$ of $M_a$ is $\mathbb{Z}_2$, but the $\pi_1$ of $M_a$ as a
   topological space is trivial.

    Theorem~\ref{thm1} is
  not true for non-compact symplectic  manifolds. For instance, take $S^1\times\mathbb{R}$
  and let $S^1$ act by rotating the first factor. This action is Hamiltonian with moment map being the projection
  to $\mathbb{R}$. We see that each reduced space is a point.\\

        In this paper, we
    generalize Theorem~\ref{thm1} to the case of torus actions and
    to the case of
    non-abelian group actions.
   We state this generalization in 2 theorems, separating the abelian and non-abelian group actions.
      \begin{theorem}\label{thm2}
     Let $(M, \omega)$ be a connected, compact symplectic manifold equipped with a Hamiltonian torus $T^n$ action ($n>1$) with moment map $\phi$.
   Then, as fundamental groups
   of topological spaces, $\pi_1(M)\cong\pi_1(M_{red})$, where $M_{red}$ is the symplectic quotient at any value  of the moment map $\phi$.
   \end{theorem}

      \begin{theorem}\label{thm3}
     Let $(M, \omega)$ be a connected, compact symplectic manifold equipped with a Hamiltonian  $G$ action with moment map $\phi$, where $G$
   is a connected compact non-abelian Lie group.
   Then, as fundamental groups
   of topological spaces, $\pi_1(M)\cong\pi_1(M_{red})$, where $M_{red}$ is the symplectic quotient at any coadjoint orbit in
   the image of the moment map $\phi$.
   \end{theorem}

   In order to explain the main ideas of the proofs, let us recall
   the following convexity theorems.

       \begin{theorem}\label{abelian-convexity}
      (\cite{AT} or \cite{GS0})
      Let $(M, \omega)$ be a connected compact symplectic manifold
      equipped with a Hamiltonian torus $T$ action. Let $\phi$ be the moment map.
      Then

      1. For each value $a\in im(\phi)$, $\phi^{-1}(a)$ is connected.

      2. The moment
      map image $\phi(M)$ is a convex polytope $\triangle$. It is
      the convex hull of the images of the fixed point sets of $T$.
     \end{theorem}

     \begin{theorem}\label{nonabelian-convexity}
      (\cite{K} or \cite{LMTW})
  Let $(M, \omega)$ be a connected, compact symplectic manifold equipped with a Hamiltonian  $G$ action
  with moment map $\phi$, where $G$
   is a connected compact non-abelian Lie group. Let $\mathfrak{t}_+^*$
   be a fixed closed positive Weyl chamber of $\mathfrak{g}^*$. Then

   1. For each coadjoint orbit $\mathcal{O}$ in the image of $\phi$, $\phi^{-1}(\mathcal{O})$ is
     connected.

   2. The image
   $\phi(M)\cap\mathfrak{t}_+^*=\triangle'$ is a convex polytope.
     \end{theorem}

  The proofs of Theorem~\ref{thm2} and Theorem~\ref{thm3} mainly consist two steps.
  Step 1,
   we use Theorem~\ref{thm1} and ``maximal"
  value on the moment polytope to prove that the reduced space at this value has the fundamental group
   of $M$.  Step 2, we prove
  that all the reduced spaces have isomorphic fundamental groups.  Step 1 is not hard.
   For  Step 2, let us mainly explain the idea for the case of an abelian
  group $G=T$ action. The moment polytope $\triangle$ in Theorem~\ref{abelian-convexity}
 consists of $\mathbf{faces}$ with different dimensions. We call the values in the
 maximal dimensional $\mathbf{faces}$ generic values. We call the values on other $\mathbf{faces}$ singular
 values.  For two generic values near
 each other, by the equivariant
   coisotropic embedding theorem, we show that the reduced spaces at
   these two values are diffeomorphic. To compare the fundamental groups of the reduced spaces
   at a singular value $c$ and at a nearby generic value $a$, we mainly use
   two facts. The first fact is, the gradient flow  of the
   moment map square gives an
   equivariant deformation retraction from $\phi^{-1}(U)$ to $\phi^{-1}(c)$, where $U$
   is a small open neighborhood of $c$ (see \cite{W}, or \cite{Ler}). Therefore,
   $\pi_1(\phi^{-1}(U)/T)\cong\pi_1(\phi^{-1}(c)/T)\cong\pi_1(M_c)$.
   The second fact is the key step toward solving the problem. It is
  a removing lemma. The space $\phi^{-1}(U)$ is a $T$-invariant smooth manifold.
  The quotient space $\phi^{-1}(U)/T$ is a stratified
  space. We will remove in a certain order the singular strata and, possibly, some piece of the generic stratum when
  necessary from this quotient. We prove that, each time we remove,
  the resulting space has isomorphic
  fundamental group as that of $\phi^{-1}(U)/T$. We do the removing until we get a space
  which has the homotopy type
  of $M_{a}$. This proves that $\pi_1(M_c)\cong\pi_1(M_a)$. For convenience, in the above argument, we may take
  a subset $\bar{U'}$ of $U$ and prove the above fact by using $\bar{U'}$ (see Lemma~\ref{lem3.4}).
  For the case of a non-abelian Lie group $G$ action, the closed positive Weyl chamber consists of
  faces with different dimensions each of which has a different stabilizer group under the
  coadjoint action.  We use the
  symplectic cross section theorem (see Theorem~\ref{thmcross}) to
  reduce the proof of Step 2 for the values on the maximal
  dimensional face of $\mathfrak{t}_+^*$ (which contains values of $\phi$) to a torus action case. For
  moment map values on other faces of $\mathfrak{t}_+^*$, we use the
   symplectic cross section theorem and a similar idea as in the
   case of an abelian group action.\\

   \begin{remark}
   By the above description (and by the method of the proof of
   Theorem~\ref{thm1} in \cite{L}), the isomorphisms of
the fundamental groups between some of the  two
   spaces are obtained by the fact that the two spaces are
   diffeomorphic, homotopy equivalent, or by using the Van-Kampen
   theorem when we do removing (or gluing). When we do removing from
   a space, we may need to do it multiple times. Each time we remove,
  we prove that the resulting space has isomorphic fundamental group
  as the previous one. Each time we use the Van-Kampen theorem, the base
  point is naturally taken in the connected intersection of the two
  connected open
  sets which cover the space. For the above reasons, in
   this paper, to simplify notation, when the context is clear,
   we will omit writing base point
   when we write $\pi_1$ of a space. The isomorphisms of the
   fundamental groups between the rest of the two spaces are
   obtained by transitivity.
    \end{remark}

   The method of the proofs of  Theorem~\ref{thm2} and Theorem~\ref{thm3}
   for the part that all the reduced spaces
   have isomorphic fundamental groups can be applied to the case of non-compact manifolds with proper moment maps.
   Regarding the fundamental group of the manifold $M$, we have the following observation.
   If the moment map $\phi$ has no critical values, then $\phi$ is a proper submersion from $M$ to $\mathfrak{g}^*$.
   By Ehresmann's Lemma, $\phi$ gives a fibration from $M$ to $\mathfrak{g}^*$ with connected fiber (\cite{K} or \cite{LMTW})
    diffeomorphic to $\phi^{-1}(a)$ for
   some $a\in im(\phi)$. By the homotopy exact sequence for fibrations, we have $\pi_1(M)\cong\pi_1(\phi^{-1}(a))$ (see the example
   of $S^1\times\mathbb{R}$ following Theorem~\ref{thm1}). This may not
   be the fundamental group of the reduced spaces. If the moment map $\phi$ has critical values, one may still able to
   prove that $\pi_1(M)\cong\pi_1(M_{red})$. There can be different approaches for this.\\

    Another interesting fact that the method  implies is:
   \begin{theorem}\label{thm4}
      Let $(M, \omega)$ be a connected, compact symplectic manifold equipped with a Hamiltonian  $G$ action
      with moment map $\phi$, where $G$
   is a connected compact Lie group.
   Then, as fundamental groups
   of topological spaces, $\pi_1(M/G)\cong\pi_1(M_{red})$, where $M_{red}$ is the symplectic quotient at any coadjoint
   orbit in  the image of the moment map $\phi$.
   \end{theorem}

   In this paper, when we say fundamental group, we mean the
     fundamental group of the topological space, without explicitly
     saying so.\\

    A brief organization of the paper: In Section 2, we will recall
    facts about proper compact Lie group actions and stratified
    spaces. Most importantly, we build blocks on removing certain strata from a stratified space
    will keep $\pi_1$ of the resulting space the same as that of the original space. In Section 3,
    we prove Theorem~\ref{thm2}. In Section 4, we recall the
    symplectic cross section theorem and the local normal form
    theorem. These are important tools for studying non-abelian Hamiltonian Lie
    group actions. In Section 5, we prove Theorem~\ref{thm3} for
    rank 1 connected compact non-abelian Lie group actions. This
    not only gives us an example of Theorem~\ref{thm3}, but also
    allows us to see the method of proof for general non-abelian
    connected compact Lie group actions. In Section 6, we prove
    Theorem~\ref{thm3}. In Section 7, we prove
    Theorem~\ref{thm4}.
   \subsubsection*{Acknowledgement} This work was supported in part by the Center for
mathematical analysis, geometry, and dynamical systems, IST,
Lisbon, Portugal where I found the method of the proof; and it
was supported in part by the research grant R1F105L15 of
professor Martin Schlichenmaier in Luxembourg University where I
wrote most part of the paper. I thank both institutes and the
people there for support and hospitality. I would like to thank
Sue Tolman for her discussion especially her explanation on
symplectic cross section theorem; and Gustavo Granja for
his nice discussion on Topology when I was studying $SU(2)$ and
$SO(3)$ actions. I thank Oleg Sheinman for some discussion on
Lie groups and for his pointing out Armstrong's theorem.
Finally, I thank the referee for his comments.

  \section{Lie group actions and stratified spaces}

      In this Section, we will recall the definition of a stratified
      space. For proper Lie group actions on a smooth manifold where slice
      theorem applies, Bredon has shown that the quotient space is a
      stratified space. We will recall this proof, and emphasize
      important points about stratified spaces. Then, we state two
      lemmas and Armstrong's theorem  which are very useful in the subsequent sections about
      removing strata from a stratified space.

       One may refer to \cite{SL} for the following definition of $\mathbf{stratified\, spaces\,}$.
       \begin{definition}\label{defdec}
   Let $X$ be a Hausdorff and paracompact topological space and let $\mathcal J$ be a partially ordered set with order relation denoted by $\leq$. A $\mathcal J$-decomposition of $X$ is a locally finite collection of disjoint, locally closed manifolds $S_i\subset X$ (one for each $i\in\mathcal J$) called pieces such that

   (i) $X=\cup_{i\in\mathcal J}S_i$;

   (ii) $S_i\cap\bar{S_j}\neq\emptyset\Leftrightarrow S_i\subset\bar{S_j}\Leftrightarrow i\leq j$.\\
  We call the space $X$ a $\mathcal J$-$\mathbf{decomposed\, space\,}$.
   \end{definition}
 \begin{definition}\label{defstr}
  A $\mathbf{decomposed\, space\,}$ $X$ is called a $\mathbf{stratified\, space\,}$ if the pieces of $X$, called strata, satisfy the following condition:

   Given a point $x$ in a piece $S$, there exist an open neighborhood $\tilde{U}$ of $x$ in $X$, an open ball $B$ around $x$ in
  $S$, a compact stratified space $L$, called the $\mathbf{link\, of\,x}$, and a homeomorphism
   $$\varphi: B\times \overset{\circ}{C}L\rightarrow \tilde{U}$$
  that preserves the decompositions.
  \end{definition}
   In the above definition, $\overset{\circ}{C}L$ is the space obtained by collapsing the boundary $L\times {0}$ of the
  half-open cylinder $L\times [0, \infty)$ to a point.

    From now on, for convenience, we may not specify the point $x$ in a (connected) stratum $S$,
    and we will call the $\mathbf{link\, of\, x}$ the $\mathbf{link\, of\, S}$.\\

   Let us recall Bredon's idea (see \cite{Br}) about the fact that the quotient space of a proper Lie group action on a smooth manifold is a stratified space.
  Assume a compact Lie group $G$ acts smoothly on a smooth manifold $N$. Assume that we have chosen a $G$-invariant metric on $N$.
  By the slice theorem, a neighborhood in $N$
 of an orbit with isotropy type $(H)$ (subgroups conjugate to $H$) is isomorphic to $G\times_H D$, where $D$ is a disk in the orthogonal complement of the tangent space
 to the orbit on which $H$ acts linearly.  The action of $H$ on $D$ is equivalent to an ``orthogonal'' action, i.e., $H$ acts as a subgroup
  of $O(n)$, where $n=$dim $D$ (See Theorem 0.3.5 in \cite{Br}). So a neighborhood of
 the point $G/H$ in $N/G$ is isomorphic to $D/H$. Let $D=D_1\times D_2$, where $D_1$ is the disk consisting
  of fixed vectors  by $H$. So the neighborhood
 of $G/H$ in $N/G$ is isomorphic to $D_1\times\overset{\circ}{C}L$, where $L=S(D_2)/H$ with $S(D_2)$ being a ``sphere'' of $D_2$. It is clear
 that the link $L$ is connected if the dimension of $D_2$ is bigger than 1, and the orbit types in $L$ have isotropy groups no bigger than $H$. So
  a small neighborhood in $N/G$ of a point in the $(H)$-stratum
 consists of $(K_i)$-strata, where $K_i$ is conjugate to a subgroup of $H$. In other words, the $(H)$-stratum can only be contained
  in the closure of $(K)$-strata, where $K$ is sub-conjugate to
  $H$. If we take a connected component $S$ of the $(H)$-stratum in
   $N/G$,
  then a neighborhood of $S$ in $N/G$ can be seen as the total space of a fiber bundle over $S$
   with fiber $\overset{\circ}{C}L$.\\

    In Section 3, Section 5 and Section 6, we will compute links of
    points. For convenience, we summarize the above computation of a
    link in the following Lemma.
    \begin{lemma}\label{lemlink}
        Let a compact connected Lie group $G$ act smoothly on a smooth
   connected manifold $N$. Then the quotient space $N/G$ is a stratified
   space. Let $S$ be a connected component of a stratum of $N/G$ with
   isotropy group conjugate to $H$. Let $G\times_H D$ be a
   neighborhood in $N$ of an orbit  with isotropy type $(H)$ ($(H)$
   denotes the collection of subgroups conjugate to $H$). Let $D=D_1\times D_2$, where $D_1$ is the disk consisting of
    fixed vectors by $H$.
   Let $S(D_2)$ be a sphere in $D_2$.
   Then the $\mathbf{link}$ of the corresponding point in $S$ (or the link of $S$) is
   $L_H=S(D_2)/H$. If dim$(D_2)>1$, then $L_H$ is connected.
   \end{lemma}

    We will mainly use the following lemma and theorem to determine whether the link of a
    connected component of a stratum is simply connected.
  \begin{lemma}\label{lembredon}
    (Corollary 6.3 in \cite{Br}) If $X$ is an arcwise connected $G$-space, $G$ compact Lie,
    and if there is an orbit which is connected, then
    the fundamental group of $X$ maps onto that of $X/G$.
    \end{lemma}
    \begin{theorem}\label{thmarm}
    (M. A. Armstrong, \cite{A})
   Let $G$ be a compact Lie group acting on a connected, locally path connected, simply connected, locally
   compact metric space $X$. Let $H$ be the smallest normal subgroup of $G$ which contains the identity component
   of $G$ and all those elements of $G$ which have fixed points.
   Then the fundamental group of the orbit space $X/G$
   is isomorphic to $G/H$.
    \end{theorem}

    Note that, if $G$ is connected, the above theorem claims the same as Lemma~\ref{lembredon} for simply connected $X$.\\

The next lemma is the key ingredient for claiming that the fundamental
group remains unchanged  after removing a
    stratum.
  \begin{lemma}\label{lemstr}
   Let a compact connected Lie group $G$ act smoothly on a smooth
   connected manifold $N$ such that $N/G$ is a stratified space.
    Let $S$ be a connected component of a stratum in $N/G$ such that no
   other strata are contained in the closure of $S$. Assume the $\mathbf{link}$
   $L_S$ of $S$ is connected and simply connected. Then
   $\pi_1(N/G)\cong\pi_1(N/G-S)$.
 \end{lemma}

   \begin{proof}
  Let $O_1$ be an open neighborhood of $S$
  in $N/G$ such that $O_1$ fibers over $S$ with fiber
 $\overset{\circ}{C}L_S$.
  Since $\overset{\circ}{C}L_S$ is simply connected, we have $\pi_1(O_1)\cong\pi_1(S)$. Take $O_2=N/G-S$.
  Then $O_1\cap O_2$ fibers over
 $S$ with fiber $\overset{\circ}{C}L_S-0$ which is homotopy equivalent to $L_S$. Since $L_S$ is connected and
 simply connected by assumption,
  $\pi_1(O_1\cap O_2)\cong\pi_1(S)$. By the Van-Kampen theorem, $\pi_1(N/G)\cong\pi_1(O_1)*_{\pi_1(O_1\cap O_2)}\pi_1(O_2)\cong\pi_1(N/G-S)$.
 \end{proof}

  \section{Proof of Theorem 2}

      By Theorem~\ref{abelian-convexity}, under the assumptions of Theorem 2,  the moment map image is a convex polytope.\\

      Now, let us recall the following well known facts about connected compact
      Lie group actions. One may refer to \cite{GS0} and the
      references cited in \cite{GS0}.
      \begin{proposition}\label{propcodm}
      Let $X$ be a connected manifold and $G$ be a connected compact
      Lie group acting smoothly on $X$. Let $G_x$ be the stabilizer
      group of $x$. Let $r=$min dim $G_x$. Let $X_i$ be the set of
      $x$'s for which dim $G_x\geq i+r$. Then

      1. If $X$ is compact, up to conjugacy, only a finite number of
      subgroups of $G$ occur as stabilizer groups of points of $X$.
       The subset $X_i$ has codimension  $\geq i+1$
      (for $i=1, 2,...$ ) in $X$.

      2. Let $H$ be a closed subgroup of $G$. Let $X_H=\{x\in X\,|\, G_x=H\}$.
      Then $X_H$ is a submanifold of $X$. If $X$ is a symplectic
      manifold and $G$ acts symplecticly, then $X_H$ is a symplectic
      submanifold of $X$.
     \end{proposition}

      From 1. of the above proposition, $X_r$ is open dense and
      connected. Up to conjugacy, let us call the stabilizer type
      of the points in $X_r$ $\mathbf{principal\, stabilizer\, type}$. \\

      Now we come back to our symplectic manifold $(M, \omega)$ with a Hamiltonian $G$ action, where $G$ is a connected
      compact Lie group. By definition of the moment map $\phi$, for each $X\in\mathfrak{g}$, we have
      $i_{X_M}\omega=d<\phi, X>$, where $X_M$ is the vector field on $M$ generated by
      $X$. From this, we easily derive

      \begin{lemma}
       Let $m\in M$. Let $G_m$ be the stabilizer group of $m$ in $G$ and let
       $\mathfrak{g}_m$ be its Lie algebra. Then the image of $d\phi_m: T_m\rightarrow\mathfrak{g}^*$ is the annihilator
       in $\mathfrak{g}^*$ of $\mathfrak{g}_m$.
      \end{lemma}

        For the manifold $M$ in Theorem~\ref{abelian-convexity},
      more explicitly, it can be stratified according to the isotropy
     groups. Let $T_1, T_2, ..., T_N$ be the subgroups of $T$ which
     occur as stabilizer groups of points of $M$. Let $M_i$ be the set
     of points in $M$ for which the stabilizer group is $T_i$. By
     relabeling, we may assume that the $M_i$'s are connected (So some $T_i$' may be repeated). Then
     $M$ is a disjoint union:
     $$M=\bigcup_{i=1}^N M_i.$$
    The moment map images of the fixed point set components of $T$
      are called the $\mathbf{vertices}$ of $\phi$. The $\mathbf{vertices}$ of $\phi$ can be
      ``real" vertices on the boundary of $\triangle$, and
      can be  ``vertices" inside $\triangle$.

      Each $M_i$ is a $T$-invariant symplectic submanifold of $M$,
      $\phi(M_i)$ is an open subset of the affine plane
      $a_i+\mathfrak{t}_i^{\bot}$, where $a_i$ is a vector in
      $\mathfrak{t}^*$, $\mathfrak{t}_i=Lie(T_i)$ and $\mathfrak{t}_i^{\bot}$ is the annihilator (or perpendicular) of
      $\mathfrak{t}_i$ in  $\mathfrak{t}^*$ (or in $\mathfrak{t}$). Moreover, $\phi(M_i)$ is the union of a finite number of convex
      sets each of which is the convex hull of a collection of the
      $\mathbf{vertices}$.

     By the above description, $\triangle$ consists of
      $\mathbf{faces}$ with different dimensions.

      \begin{remark}\label{remeffective}
       We may assume that $\triangle$ contains an open subset of $\mathfrak{t}^*$.  By the definition of
       the moment map, this is
       the same as assuming that the $T$ action has finite generic stabilizer group. Moreover, we assume that
      $\bigcap_{m\in M} T_m=1$, where $T_m$ is the stabilizer group of $m$. If this is not the case, we  divide $T$ by the
       common stabilizer group and consider the quotient torus action.
      Let us call the values of $\phi$ in the
       open set $\mathbf{regular\, values}$.
     Clearly, the set of $\mathbf{regular\,
       values}$ is open and dense in  $im(\phi)$.
      \end{remark}

      Let us call a connected set of
      regular values of $\phi$ a (connected) $\mathbf{chamber}$ of
      $\triangle$.  The moment map image  $\triangle$ may have one
      or more than one connected $\mathbf{chambers}$.

     \begin{remark}\label{rem1}
     Note that the $M_i$'s are disjoint, but the $\phi(M_i)$'s may
     not  be disjoint (they are disjoint if $(M, \omega, T, \phi)$ is a  completely integrable system.).
     For instance, regular points on $M$, i.e., points with finite stabilizer groups
     can be mapped to non-open $\mathbf{faces}$. One may see this
     easily for $S^1$ actions, and then generalize to $T$ actions.
     For a boundary $\mathbf{vertex}$ $v$, $\phi^{-1}(v)$ only consists of one
     fixed point set component; for a $\mathbf{vertex}$ inside $\triangle$,
     this cannot be true. Indeed, by the following Theorem~\ref{thmform}, a neighborhood in $M$ of a fixed point $x$
     is isomorphic to a $T$ representation $V=W\oplus V^T$, and $\phi|_V$ is the moment map for the $T$ action on $W$.
     If $\phi(x)$ is on the boundary of $\triangle$, then at least one subcircle of $T$ acts  on $W$ with  weights
     all positive, so  $\phi^{-1}(0)$ only consists of the fixed point set component containing $x$. One may see similarly
     that for an interior vertex $v$, $\phi^{-1}(v)$ contains more points than the fixed point set itself.
     For a $\mathbf{wall}$ $\mathcal W$ on the boundary,
     $\phi^{-1}(\mathcal W)$ is fixed by the circle generated by the
     direction perpendicular to $\mathcal W$. For a $\mathbf{wall}$ $\mathcal W$
      inside $\triangle$,   $\phi^{-1}(\mathcal W)$ consists more than just a submanifold which is fixed by the circle generated by
      the direction perpendicular to $\mathcal W$.
     \end{remark}

   \begin{lemma}\label{lem3.1}
    Suppose $X_1, X_2,..., X_k$ are connected components of the
    fixed point set of $T$ such that $\phi(X_i), i=1,...,k$ are
    boundary vertices $v_i, i=1,..., k$ of the polytope $\triangle$. Then
    $\pi_1(M_{v_i})\cong\pi_1(X_i)\cong\pi_1(M)$, for $i=1,..., k$.
   \end{lemma}
   \begin{proof}
   The first equality is clear since $\phi^{-1}(v_i)=X_i$ (see Remark~\ref{rem1}) and $X_i$ is fixed by $T$. To prove the second equality, choose
   a sub-circle in $T$ such that $v_i$ is the maximal value of the
   moment map for the circle action. Then we apply
   Theorem~\ref{thm1}.
   \end{proof}

        \begin{remark}\label{circle}
    Although we are proving Theorem~\ref{thm2}, the following proofs of Lemma~\ref{lem3.2}
    and Lemma~\ref{lem3.3} cover the case of $S^1$ actions.
     \end{remark}

   \begin{lemma}\label{lem3.2}
      For two values $a, b$ near each other in one connected $\mathbf{chamber}$ of $\triangle$, we have
     $\pi_1(M_a)\cong\pi_1(M_b)$. Therefore, by connectivity of the $\mathbf{chamber}$, for all values $a$ in
     this $\mathbf{chamber}$,
    $\pi_1(M_a)'$s are all isomorphic.
    \end{lemma}

    \begin{proof}
     Take two regular values $a$
     and $b$ close enough. By the equivariant co-isotropy embedding
     theorem, there exists a small neighborhood $U$ containing $a$ and $b$ such that
     $U$ consists of regular values, and
     $\phi^{-1}(U)$ is isomorphic to $\phi^{-1}(a)\times U$, where $T$ acts on $\phi^{-1}(a)$ and the  moment map is the
     projection to $U$. So $M_a$ is diffeomorphic to $M_b$. So $\pi_1(M_a)\cong\pi_1(M_b)$.
    \end{proof}

    \begin{lemma}\label{lem3.3}
    Let $c$ be a non-regular
    value on $\triangle$. Let $a$ be a regular value such that it is very near
    $c$. Then $\pi_1(M_c)\cong\pi_1(M_a)$.
   \end{lemma}

    Proof of Theorem~\ref{thm2}:
    \begin{proof}
   Theorem~\ref{thm2} follows from Lemma~\ref{lem3.1},
   Lemma~\ref{lem3.2}, and Lemma~\ref{lem3.3}.
   \end{proof}

     Lemma~\ref{lem3.3} follows from the following Lemma~\ref{lem3.4}.
  Let us first recall the following theorem on the convergence of the gradient flow of the moment map square:
  \begin{theorem}\label{thmLW} (\cite{W} or \cite{Ler})
           Let $(M, \omega)$ be a connected  Hamiltonian $G$-manifold with proper moment map $\phi$,
           where $G$ is a connected compact Lie
           group. Choose a $G$-invariant metric on $M$. Assume that the moment map image intersects a neighborhood of
           $0$ (the image not necessarily fills an open neighborhood of
           $0$). Then there exists a $G$-invariant open neighborhood $U\subset\mathfrak g^*$  of
           $0$ such that the negative gradient flow of the moment map square
           induces a $G$-equivariant deformation retraction
          from  $\phi^{-1}(U)$ to $\phi^{-1}(0)$.
  \end{theorem}

       If $G$ is a torus, we can always shift the moment map (by a constant) such
       that a value $c$ we consider corresponds to the $0$ value of
       the new moment map, so without loss of generality, we can
       regard $c$ as $0$.

       Assume $G$ is a torus. Take $U$ as in Theorem~\ref{thmLW}.
       Let
       $U'$ be the intersection of $U$ with a connected open $\mathbf{chamber}$, and,
       let
       $\bar{U'}$ be its closure in $U$. Then $\phi^{-1}(\bar{U'})$ also $G$-equivariantly deformation
       retracts to $\phi^{-1}(0)$.

     \begin{lemma}\label{lem3.4}
      Let $c$ be a value on a singular $\mathbf{face}$ of $\triangle$, and let $a$ be a regular value very near $c$.
      Let $U$ be a small open neighborhood of $c$ on $\triangle$ containing $a$ such that $\phi^{-1}(U)$ equivariantly
      deformation retracts to
      $\phi^{-1}(c)$. Let $U'$ be the intersection of $U$ with the connected open $\mathbf{chamber}$
      containing $a$, and let
     $\bar{U'}$ be its closure in $U$.
      Let $B$ be the set of values in $\bar{U'}$ but not in $U'$.
      Then
      $\pi_1(\phi^{-1}(\bar{U'})/T)\cong\pi_1(\phi^{-1}(\bar{U'})/T-\phi^{-1}(B)/T)$, i.e.,
      $\pi_1(M_c)\cong\pi_1(M_a)$.
     \end{lemma}

      We could use the set $U$ itself, and apply a removing and flowing (by using the gradient flow) process to
      achieve $\pi_1(M_c)\cong\pi_1(M_a)$. I found that I still would have to do the above removing in the end. So
      taking  $\bar{U'}$ is more convenient. \\

     In order to prove the above lemma, let us recall the following Local Normal Form theorem for abelian Lie
     group actions.

         \begin{theorem}\label{thmform}
        (Local normal form) (\cite{GS1}) Let $(M, \omega)$ be a symplectic manifold with a Hamiltonian
        torus $T$ action. Let
         $H$ be the isotropy subgroup of a point $p$ in $M$. Then a neighborhood in $M$ of the orbit through $p$
         is equivariantly symplectomorphic to
         $T\times_H(\mathfrak{b}^{\bot}\times V)$, where $\mathfrak{b}^{\bot}$ is the annihilator
         of $\mathfrak{b}=$Lie$(H)$ in $\mathfrak t^*$ on which $H$ acts by the coadjoint action (trivial in this case),
     and $V$ is a complex vector
         space on which $H$ acts linearly and symplectically.

      The equivalence relation on  $T\times_H(\mathfrak{b}^{\bot}\times V)$ is given by
      $(t, a, v)\approx (th^{-1}, a, h\cdot v)$ for $h\in H$.

The $T$ action on this local model is
  $t_1\cdot [t, a, v]=[t_1t,a,v]$, and the moment map on this local model is
       $\phi([t, a, v])=\phi(p)+a+\psi(v)$, where $\psi(v)$ is the moment map for the $H$ action on $V$.
       \end{theorem}

      \begin{remark}\label{order}
      By Theorem~\ref{thmform}, if an orbit has stabilizer $H$, then the nearby orbits of this orbit have
      stabilizers no bigger than $H$. So, when we
      remove a connected  stratum which is more singular, we will not destroy the link of its nearby (less singular) strata.
      If a face $\mathcal F$ of $\triangle$ is in the closure
      of the face $\mathcal F'$, then $\phi^{-1}(\mathcal F)$
      contains more singular (with bigger stabilizer groups) strata
      than the strata in $\phi^{-1}(\mathcal F')$.  In Lemma~\ref{lem3.4},
      if $\mathcal F$ is the $\mathbf{face}$ in $U$ containing the singular value $c$, then it is the most degenerate face
      in $U$.
     \end{remark}

     The proof of Lemma~\ref{lem3.4} is a removing process. According to Lemma~\ref{lemstr}, the quotient space of a smooth manifold
     by a compact Lie group action is a stratified space, and, certain removing of strata
     from the quotient keeps $\pi_1$ of the resulting space the same as $\pi_1$ of this quotient itself. The $\phi^{-1}(\bar{U'})$ we took
     is not a manifold. Nevertheless, the analysis of the neighborhoods allows us to perform the removing in
     $\phi^{-1}(\bar{U'})/T$.

     \begin{proof}
      Assume that  $\mathcal F$ is the singular face in $U$ containing $c$.  Let $m\geq 0$ be its dimension.
      The face $\mathcal F$ is the most degenerate face in $U$ and in $\bar{U'}$. According to Remark~\ref{order}, we first
      remove $\phi^{-1}(\mathcal F)/T$.

      Assume the dimension of the torus $T$ is $n$.
      Let us identify the Lie algebra $\mathfrak{t}$ of $T$ with its
      dual $\mathfrak{t}^*$. The $m$-dimensional linear subspace $L$
      which contains  $\mathcal F$ as an open set generates a
      subtorus $T^m$. If $\mathcal F$ is not a $\mathbf{wall}$, the
complementary linear subspace of $L$
      in
      $\mathfrak{t}$ spanned by the set $S$ of other directions of
      one-dimensional $\mathbf{faces}$ of $\bar{U'}$ generates a
      complementary subtorus $T^{n-m}$.
      If $\mathcal F$ is a $\mathbf{wall}$, then we take the
      direction orthogonal to $L$ as a complementary direction, and
      we take the complementary subtorus of $T^m$ generated by this
      direction. So, we have chosen a splitting $T=T^{n-m}\times T^m$.

      Notice that the most singular stratum in
      $\phi^{-1}(\mathcal F)$ has stabilizer  $T^{n-m}\times \Gamma$, where $\Gamma$ is a
     finite subgroup of $T^m$. The set
      $\phi^{-1}(\mathcal F)$ may contain other strata with
      stabilizer(s) of the form $H=(T_1\times\Gamma')\times\Gamma$, where $T_1$ is a connected subgroup
     of  $T^{n-m}$ generated by directions of some 1-dimensional $\mathbf{faces}$ of $\bar{U'}$,
     $\Gamma'$ is a finite subgroup of a complementary group $T_2$
     of $T_1$ in $T^{n-m}$. (We may have a different finite group $\Gamma$ which is a subgroup
     of the previous $\Gamma$. But this will not affect the proof, so we use the same notation). Here, we chose $T_2$ to be the
     subgroup of $T^{n-m}$ generated by the rest of the directions in $S$ of
     1-dimensional $\mathbf{faces}$ of $\bar{U'}$ when $\mathcal F$ is not a $\mathbf{wall}$,
     and, we chose $T_2$ to be trivial or to be the subgroup generated
     by the orthogonal direction of $L$ when  $\mathcal F$ is a
     $\mathbf{wall}$ (remember that in this case $n-m$ was $1$).

              By Theorem~\ref{thmform}, a
        neighborhood in $M$ of an orbit with
        isotropy type $H$ is isomorphic to
        $A=T\times_H(\mathbb{R}^l\times\mathbb{R}^m\times V)$. The moment map
        on $A$ is $\phi=a+b+\psi(v)$, where $a\in \mathbb{R}^l$, $b\in
        \mathbb{R}^m$, and $\psi$ is the moment map of the $H$
        action  on $V$
         (we assumed that $``c=0"$).
        We split $V=W\times V^H$, where $W$ has no non-zero
    fixed vectors by $H$. Then $\psi$ is just $\psi|_W$. So
    $A\cap\phi^{-1}(\bar{U'})=T\times_H((\mathbb{R}^m\times V^H)\times((\mathbb{R}^+)^l\times W\cap\psi^{-1}(\bar
     {U'})))$ (strictly speaking, the above $\mathbb{R}^m$ should be
     a small open disk in
     $\mathbb{R}^m$ corresponding to $\mathcal F\cap\bar{U'}$),
     where $(\mathbb{R}^+)^l$ are the non-negative real (half) lines
    pointing towards $\bar{U'}$.
    The $H$-stratum in $A\cap\phi^{-1}(\bar{U'})$ which was mapped to
     $\mathcal F$ is
        $T\times_H(\mathbb{R}^m\times V^H)$.
     The link $L_H$ of the corresponding quotient $H$-stratum  in $(A\cap\phi^{-1}(\bar{U'}))/T$ is
     $S((\mathbb{R}^+)^l\times W\cap\psi^{-1}(\bar{U'}))/H$. Here,
        $S((\mathbb{R}^+)^l\times W\cap\psi^{-1}(\bar{U'}))$ is the
    intersection of  $S(\mathbb{R}^l\times W)$ with $(\mathbb{R}^+)^l\times W\cap\psi^{-1}(\bar{U'})$.
     Now, we consider all the possible cases of $H$ (will
     correspond to different $l$).

       1. In the case of $l=0$ (corresponding to $H=T^{n-m}\times\Gamma$),
       since we assumed that the moment map value fills  $U'$,
    the moment map $\psi|_W$ has to be
       non-trivial.
       By Remark~\ref{finite} and Lemma~\ref{ints} below, the
       link is connected and simply connected.

       2. In the case of $l\neq 0$ and $W\neq 0$, by Remark~\ref{rem ints'} and
       Lemma~\ref{ints'} below, the link $L_H$ is
     connected and simply connected.

       3. In the case of $l\neq 0$ and $W=0$ (all
    the points in $A\cap\phi^{-1}(\bar{U'})$ have the same stabilizer group
    $H$ which has to be the generic stabilizer group), the link
    $L_H=S((\mathbb{R}^+)^l)/H=S((\mathbb{R}^+)^l)$
    is connected and simply connected.

      If there are more singular faces left in $\bar{U'}$, we remove similarly as above.

     \end{proof}

      \begin{remark}\label{finite}
     For Case 1 in the proof of Lemma~\ref{lem3.4}, we needed to
     consider the quotient
     $(S(W)\cap\psi^{-1}(\bar{U'}))/(T^{n-m}\times\Gamma)$,
     where $W$ is a complex $(T^{n-m}\times\Gamma)$-representation
     isomorphic to some $(\mathbb{C})^n$ (it splits into a product of
     $\mathbb{C}$) on which $T^{n-m}\times\Gamma$ acts as a subgroup
     of the maximal torus of $U(n)$.
     Due to how cyclic group acts on $(\mathbb{C})^n$, and due to
     the fact that  a finite group action does not contribute
     to $\psi$, we can first divide $W$ by $\Gamma$, we get
     $W/\Gamma$ homeomorphic to $W$. The action of $T^{n-m}$ on  $W/\Gamma$
     corresponds to a ``weight change" comparing to the action of
     $T^{n-m}$ on $W$. Therefore, we can restrict attention to the case of the
     following lemma.
       \end{remark}

   \begin{lemma}\label{ints}
         Let $\mathbb{C}^n$ be an effective $T^k$ symplectic representation, where $T^k$ is a connected torus, and $k\leq n$.
     Let $\psi$ be the moment map for the
     $T^k$ action.  Let $U'\subset (\mathfrak{t}^k)^*$ be an open connected
     $\mathbf{chamber}$ consisting of regular values of $\psi$, and let $\bar{U'}$ be its closure.
     Let $S'=S^{2n-1}\cap \psi^{-1}(\bar{U'})$.
     Then, the quotient $S'/T^k$
      is always connected and simply connected.
     \end{lemma}
     \begin{proof}
     When $k=n$, $T^n$ corresponds to the maximal torus of $U(n)$. It acts on
      $\mathbb{C}^n$ in the standard way with
      moment map $\psi=(|z_1|^2, |z_2|^2, ... , |z_n|^2)$. In this case,  $S'=S^{2n-1}$. So $S'/T^k$ is
      either a point when $n=1$, or it is connected and simply connected by Lemma~\ref{lembredon}.

      Now assume that $k<n$. We may assume that $n>1$.
      Then $T^k$ acts on $\mathbb{C}^n$ as a subtorus of $T^n$.
     Its moment map $\psi$ is the projection of the above moment map to the dual Lie algebra of $T^k$, i.e.,
     $\psi=\alpha_1|z_1|^2+\alpha_2|z_2|^2+...+\alpha_n|z_n|^2$, where $\alpha_i, i=1,..., n$ are weight vectors
     in $(\mathfrak{t}^k)^*$. Then  $\bar{U'}$  is formed by the cone with non-negative coefficients spanned by
     $p$ with $p\geq k$ number of  vectors among $\alpha_i, i=1,..., n$. Without loss of generality, we assume that they are the
     first $p$ vectors. Any $k$ number of linearly independent vectors among them
     generates $T^k$. Take $k$ number of linearly independent vectors, say the first $k$ vectors, among the $p$ number of vectors, and
       write each of $\alpha_i, i=p+1,..., n$ as linear combinations of $\alpha_i, i=1,..., k$:
      $\alpha_i=a_{i1}\alpha_1+...+a_{ik}\alpha_k$ for $i=p+1,..., n$. Then
      $\psi=(|z_1|^2+\sum_{i\geq p+1} a_{i1}|z_i|^2)\alpha_1+...+(|z_k|^2+\sum_{i\geq p+1} a_{ik}|z_i|^2)\alpha_k+|z_{k+1}|^2\alpha_{k+1}+...+|z_p|^2\alpha_p
      =\sum_{i\leq p}f_i\alpha_i$.
      So $\psi^{-1}(\bar{U'})=\{z\in\mathbb{C}^n: f_i(z)\geq 0, i=1,...p\}$. Since the action is linear and the moment map
      is homogeneous, we only need to prove that $((\mathbb{C}^n-0)\cap\psi^{-1}(\bar{U'}))/T^k$ is connected and simply
      connected. Now consider $(\mathbb{C}^n-0)\cap\psi^{-1}(\bar{U'})$. Since $0$ (very singular)
      is taken away,
      we can perturb the set a little bit using the gradient flow of some $f_i$ without changing the topology so
      that  the above intersection has the same topology
      as a union of the sets $A_i=\{z\in\mathbb{C}^n: f_i>0,f_J\geq 0\}$ for some $i=1, ..., p$, where
       $J=\{1,..., i-1, i+1,..., p\}$.
      We move the terms with negative coefficients in each of $f_i$ to
      the right hand side of $f_i\geq 0$. For a fixed set $A_i$, we
       equivariantly deformation retract all the $z_ls'$ (which occured in $f_i>0$ and in $f_J\geq 0$)
      with $l>p$ which are not
       on the left hand side of $f_i$ to $0$. Now, the possibilities
       are:
      (a). The right hand sides of all
       $f_J$ are $0$, then the inequalities $f_J\geq 0$ do not give
       any condition. The inequality $f_i>0$ (after deforming the right hand side to $0$)
       represents a copy of $\mathbb{C}^*$ or
       a simply connected sphere. So $A_i$ is deformed into $\mathbb{C}^*$ or
       a simply connected sphere or a product of one of them with
       some copies of $\mathbb{C}$ represented by some coordinates
       on the left hand sides of $f_J$ but not on the left hand side of $f_i$.
       So $A_i/T^k$ is connected and simply
       connected.
      (b). The right hand sides of $f_J$ are not $0$. Then, we
      perform some algebraic operations between the inequalities,
      we have a new $f_i$. We equivariantly deform the coordinates
      which are not on the left hand side of the new $f_i$ to $0$
      again, and, we may need to repeat this process until we have case (a) with a different $f_i$ from
      the original one.
    The set $((\mathbb{C}^n-0)\cap\psi^{-1}(\bar{U'}))/T^k$ is obtained
   by gluing the connected and simply
      connected $A_i/T^k$s'. Now, we only need to see that the intersection
      of each two of $A_is'$ is connected. We may treat the intersections similarly
      as the above (now we have two strict inequalities), i.e., we use algebraic operations
       and we use
      deformations until we can deform all the right hand sides to
      $0$. Now we have two strict inequalities $f_i>0, f_j>0$ with
      (quadratic terms and)
      positive coefficients. The  inequalities represented by $f_{J'}\geq 0$ (where $J'=\{1, ..., p\}-\{i, j\}$)
     with positive coefficients do not give conditions.
      So the coordinates on the left hand side of $f_{J'}\geq 0$ but not on the left hand sides of $f_i>0, f_j>0$ are free.
      If there are no common coordinates in $f_i>0$ and $f_j>0$, we see that the intersection
      is a product of connected sets therefore connected. Otherwise,
      by writing the intersection as a union of products of $\mathbb{C}^*s'$ with $\mathbb{C}s'$
      which are connected and which have connected intersections, we see that,
      the intersection of the two $A_is$ are connected. The Van-Kampen theorem justifies the conclusion.
      \end{proof}

      \begin{remark}\label{rem ints'}
      In the local model $A=T\times_H(\mathbb{R}^l\times\mathbb{R}^m\times
      V)$ with $H=T_1\times\Gamma'\times\Gamma$
      we considered in the proof of Lemma~\ref{lem3.4}, if
      we divide $A$ by $T^m$, we get
      $A/T^m=T^{n-m}\times_{(T_1\times\Gamma')}((\mathbb{R}^l\times
      W/\Gamma)\times (\mathbb{R}^m\times V^H))=(T^{n-m}\times_{(T_1\times\Gamma')}(\mathbb{R}^l\times
      W/\Gamma))\times (\mathbb{R}^m\times V^H)$. Now, for Case 2 in
      the proof of Lemma~\ref{lem3.4},
      we ``forget" the component $(\mathbb{R}^m\times V^H)$ (it
      was mapped to the $\mathbf{face}$ $\mathcal F$)
      and we count the fact that $W/\Gamma$ is
      homeomorphic to $W$ and the fact that a finite group does not contribute to the moment map,
      we may restrict attention to the case of the following lemma.
      In the following lemma, the $U'$ is the intersection of the
      $U'$ in Lemma~\ref{lem3.4} and the moment map image on $A'$.
      \end{remark}

      \begin{lemma}\label{ints'}
      Let $T$ be a $d$-dimensional torus acting in a local model
   (see Theorem~\ref{thmform})
      $A'=T\times_{(T_1\times\Gamma')}(\mathbb{R}^l\times W)$, where
      $T_1$ is a connected subgroup with dimension less than $d$ (so $l\neq 0$),
       $\Gamma'$ is a finite subgroup, and $W$ is a complex $T_1\times\Gamma'$
      representation with no non-zero fixed vectors.  Let the moment map of
      the $T$ action on $A'$ be $\phi([t, a, z])=a+\psi(z)$. Let $U'\subset \mathfrak{t}^*$ be an open
      connected
     $\mathbf{chamber}$ consisting of regular values of $\phi$, and let $\bar{U'}$ be its closure.
     Let $S'=S((\mathbb{R}^+)^l\times W\cap\psi^{-1}(\bar{U'}))$.
     Then, the quotient $S'/(T_1\times\Gamma')$
      is connected and simply connected.
     \end{lemma}
     \begin{proof}
     We can see (as we did before) that  $S'/(T_1\times\Gamma')$
      is the link of the quotient $(T_1\times\Gamma')$-stratum in
      $(A'\cap\phi^{-1}(\bar{U'}))/T$.

       The group $T_1\times\Gamma'$ acts on $(\mathbb{R}^+)^l$
       trivially; and, for a similar reason as we made in Remark~\ref{finite}
        about cyclic group action on  $W$ and about its trivial contribution
        to $\psi$, we can ``ignore" $\Gamma'$, and consider
        $S'/T_1$. Assume $T_1$ is of dimension $k$ ($k+l=d$). We use the same
        notations as we did in the proof of Lemma~\ref{ints}, we
        have (assume $W$ is isomorphic to some $\mathbb{C}^n$)
       $W\cap\psi^{-1}(\bar{U'})=\{z\in\mathbb{C}^n: f_i(z)\geq 0, i=1,...p\}$.
       For a similar reason as in the proof of
       Lemma~\ref{ints}, we consider the intersection
       $S''=((\mathbb{R}^l\times\mathbb{C}^n)-0)\cap((\mathbb{R}^+)^l\times\{z\in\mathbb{C}^n: f_i(z)\geq 0,
       i=1,...p\})$, and we prove that $S''/T_1$ is connected and simply
       connected. We perturb the set $S''$ such that it has the same
       topology as a union of the sets $A_{i',i}=\{(a_1, ..., a_l)\in\mathbb{R}^l, z\in\mathbb{C}^n: a_{i'}>0, a_{J'}\geq 0,
       f_i>0,f_J\geq 0\}$ for some $i'=1, ..., l$ or for some $i=1, ...,
       p$,
       where $J'=\{1,..., l\}-i'$ or $J=\{1, ..., p\}-i$ (only one of $i'$ and $i$ is non-empty, so
       $J'=\{1,..., l\}$ or $J=\{1, ..., p\}$).
       We argue as we did in the proof
       of Lemma~\ref{ints} that each $A_{i',i}/T_1$ is connected and
       simply connected, and the intersection of each two of these sets is connected.
        We glue them together and we use the Van-Kampen theorem to prove that $S''/T_1$
       is connected and simply connected.
       \end{proof}

  \section{Cross Section theorem and local normal form theorem}
   In this section, we will first recall the Cross-Section Theorem due to
   Guillemin and Sternberg. These cross sections will give us
   symplectic submanifolds with lower dimensional subgroup actions.
   Then, we will state the Local Normal Form Theorem for Hamiltonian Lie group actions, due to
   Guillemin-Sternberg, and Marle.

  \subsection{Cross Section Theorem}

  \begin{definition}\label{def 4.1}
   Suppose that a group $G$ acts on a manifold $M$. Given a point $m$ in $M$ with isotropy group $G_m$, a submanifold $U\subset M$ containing $m$ is a
   $\mathbf{slice\, at\, m}$ if $U$ is $G_m$-invariant, $G\cdot U$ is a neighborhood of $m$, and the map

          $G\times_{G_m}U\rightarrow G\cdot U$, \qquad      $[a, u]\longmapsto a\cdot u$   is an isomorphism.
   \end{definition}

   For instance, consider the co-adjoint action of $G=SU(2)$ or $SO(3)$ on $\mathbb{R}^3=$Lie$(G)$. Recall that all the co-adjoint orbits through non-zero points
   in $\mathbb{R}^3$ are diffeomorphic to $S^2$ (these are generic co-adjoint orbits), and the co-adjoint orbit through $0$ is $0$ (this is a singular
   co-adjoint orbit).
   For $x\in \mathbb{R}^3, x\neq 0$, there is a unique ray $I_x$
   passing through $0$ and $x$. It is easy to see that the open ray $I_x^{\circ}=I_x-0$ is a slice at $x$.
   If $x=0$, then a slice at $0$ is  $\mathbb{R}^3$.

    More generally, let us consider the co-adjoint action of a connected compact
  Lie group $G$ on $\mathfrak g^*$. For $x\in \mathfrak g^*$, let $U_x$ be the natural slice at $x$ for the
 co-adjoint action. Fix a (closed) positive Weyl chamber $\mathfrak t_+^*$, without loss of generality, we assume
  $x\in\mathfrak t_+^*$. Let $\tau\subset\mathfrak t_+^*$ be the open face of  $\mathfrak t_+^*$ containing
 $x$  and let $G_x$ be the isotropy group of $x$ (all the points on $\tau$ have the same isotropy group). Then
$U_x=G_x\cdot\{y\in\mathfrak t_+^* |G_y\subset G_x\}=G_x\cdot\bigcup_{\tau\subset\bar{\tau'}}\tau'$,
   and it is an open subset of $\mathfrak g_{\tau}^*=\mathfrak g_x^*$.\\

    We have the following Cross-Section Theorem due to Guillemin and Sternberg
(Theorem 26.7 in \cite{GS}; for the following version, see Corollary
  2.3.6 in \cite{GLS}).

    \begin{theorem}\label{thmcross}
    (Cross-Section). Let $(M, \omega)$ be a symplectic manifold with a moment map $\phi: M\rightarrow \mathfrak{g}^*$ arising from an action of a compact
    connected Lie group $G$. Let $x$ be a point in $\mathfrak{g}^*$ and let $U$ be the natural slice at $x$. Then the cross-section $R=\phi^{-1}(U)$ is
     a $G_x$-invariant symplectic submanifold of $M$, where $G_x$ is the isotropy group of $x$. Furthermore, the restriction $\phi|_R$ is a moment map
     for the action of $G_x$ on $R$.
     \end{theorem}

 \subsection{Local Normal Form Theorem}
    The following Local Normal Form Theorem describes up to
    equivariant isomorphism a neighborhood of an isotropic orbit in
    a Hamiltonian $G$-manifold.

       \begin{theorem}\label{thmform'}
        (Local normal form) (\cite{GS1}) Let $(M, \omega)$ be a symplectic manifold with a Hamiltonian connected compact Lie group $G$ action. Assume $p\in M$,
      and the orbit $G\cdot p$ is isotropic. Let
         $H$ be the isotropy subgroup of $p$. Then a neighborhood of the orbit through $p$ in $M$ is equivariantly symplectomorphic to
         $G\times_H(\mathfrak{b}^{\circ}\times V)$, where $\mathfrak{b}^{\circ}$ is the annihilator of $\mathfrak{b}=$Lie$(H)$
    in $\mathfrak g^*$ on which $H$ acts by the coadjoint action, and $V$ is
 a complex vector
space on which $H$ acts linearly and symplectically.

      The equivalence relation on  $G\times_H(\mathfrak{b}^{\circ}\times V)$ is given by $(g, a, v)\approx (gh^{-1}, h\cdot a, h\cdot v)$ for $h\in H$.

 The $G$ action on this local model is $g_1\cdot [g, a, v]=[g_1g, a, v]$, and the moment map on this local model is
       $\phi([g, a, v])=Ad^*(g)(\phi(p)+a+\psi(v))$, where $\psi(v)$ is the moment map for the $H$ action on $V$.
       \end{theorem}

\section{Proof of Theorem 3 in the case of $G=SU(2)$ or $G=SO(3)$}

   As an example, we will prove Theorem 3 for the action of $G=SU(2)$
  and $G=SO(3)$. This will give us the flavor of the proof of Theorem 3.\\

    If $0$ is the only value in the moment map image, then by the
    definition of the moment map, $G$ acts trivially on $M$. So the
    theorem is trivial in this case.\\

     Let us now assume that $\phi$ has non-zero values. Using
     Theorem~\ref{thmcross}, we will first reduce the proof of the
     theorem at non-zero values to a circle action case.

           \begin{lemma}\label{lem5.0}
        Let $(M, \omega)$ be a connected, compact symplectic manifold equipped with a Hamiltonian $SU(2)$ or $SO(3)$ action. Assume the
 moment map image has non-zero values. Then,
 as fundamental groups
   of topological spaces, $\pi_1(M)\cong\pi_1(M_{red})$,
 where $M_{red}$ is the symplectic quotient at any non-zero value  of the moment map $\phi$.
   \end{lemma}

   \begin{proof}
     Let $x\in im(\phi)$, and $x\neq 0$. Let $I_x^{\circ}$ be the natural slice at $x$,
     and take $R=\phi^{-1}(I_x^{\circ})$. By Theorem~\ref{thmcross},
     $R$ is a symplectic submanifold with a Hamiltonian $S^1$ action
     whose moment map is  $\phi|_R$ and $\phi|_R$ is  proper onto its image.
     From Theorem~\ref{nonabelian-convexity}, we can deduce that $R$ is
     connected. By Lemma~\ref{lem3.2}
    and Lemma~\ref{lem3.3} (see Remark~\ref{circle}), $\pi_1(M_a)\cong\pi_1(M_b)$ for
     any $a, b\in I_x^{\circ}$.

      This same $S^1$  acts on $M$ with moment map  being the projection of $\phi$ to $\mathbb{R}=$Lie$(S^1)$.
      Let us use $\phi_p$ to denote this
 ``projected''  moment map. Suppose it has maximal value $z$.
     Then $z\in I_x^{\circ}$, and  $\phi_p^{-1}(z)=\phi|_R^{-1}(z)$ is the maximum on $M$ of $\phi_p$.
      By Theorem~\ref{thm1},
     $\pi_1(M)\cong\pi_1(\phi_p^{-1}(z))\cong\pi_1(\phi_p^{-1}(z)/S^1)\cong\pi_1(\phi|_R^{-1}(z)/S^1)\cong\pi_1(M_z)$.

     \end{proof}

              \begin{lemma}\label{lem5.3}
     Let $(M, \omega)$ be a connected, compact symplectic manifold equipped with a Hamiltonian $SU(2)$ or $SO(3)$ action.
     Assume $0$ and $\epsilon$ small are in the moment map image. Then $\pi_1(M_0)\cong\pi_1(M_{G\cdot\epsilon})$.
          \end{lemma}

  The proof of Lemma~\ref{lem5.3} relies on Theorem~\ref{thmLW} and  Lemma~\ref{lemremove5} below.

  \begin{lemma}\label{lemremove5}
       Under the assumptions of Lemma~\ref{lem5.3},
           there exists a small $G$-invariant open neighborhood $U=\{\,x\in\mathbb{R}^3|
           \,|x|<\epsilon_0\}$, such that $U-0$ consists of regular
           values, and such that
       $\pi_1(\phi^{-1}(U)/G)\cong\pi_1(\phi^{-1}(U)/G-M_0)$.
   \end{lemma}

   Assume we have this. The proof of  Lemma~\ref{lem5.3} goes as the following.
         \begin{proof}
      Assume that we have taken $U$ small enough such that we can use Theorem~\ref{thmLW}.
      Therefore, $\pi_1(M_0)\cong\pi_1(\phi^{-1}(U)/G)$.
       By Lemma~\ref{lemremove5},  $\pi_1(\phi^{-1}(U)/G)\cong\pi_1(\phi^{-1}(U)/G-M_0)$.
      The space $\phi^{-1}(U)/G-M_0$ is homotopy equivalent
           to $\phi^{-1}(G\cdot\epsilon)/G=M_{\epsilon}$. To see this, by the Symplectic Cross Section Theorem,
            $\phi^{-1}(U)-\phi^{-1}(0)$ is the
 total space of a fibration over $S^2$ (a coadjoint orbit) with fiber $\phi^{-1}(I)$, where $I$ is an open interval. Since all values in $I$ are regular,
  $\phi^{-1}(I)$ is isomorphic to $\phi^{-1}(\epsilon)\times I$ which is equivariantly
 homotopy equivalent to $\phi^{-1}(\epsilon)$. Therefore, $\phi^{-1}(U)-\phi^{-1}(0)$ is equivariantly
 homotopy equivalent to the total space of a fibration over $S^2$ with fiber $\phi^{-1}(\epsilon)$, and this space is $\phi^{-1}(G\cdot\epsilon)$.
  So $\pi_1(M_0)\cong\pi_1(M_{\epsilon})$.
           \end{proof}

It remains to prove Lemma~\ref{lemremove5}.
\begin{proof}
    By equivariance of the moment map and by
   continuity, we see that if $\phi$ takes value $0$ and a non-zero value, it has
   to take values in an open neighborhood of $0$. Since $M$ is
   compact, by considering a subcircle action, we see that
    if we take $U$  small enough, then $U-0$
   consists of regular values of $\phi$.

    By Lemma~\ref{lemlink}, $\phi^{-1}(U)/G$ is a stratified space.
    By \cite{SL}, $\phi^{-1}(0)$ and $M_0$ are stratified spaces.
    If there is only one stratum, then these spaces are manifolds.  The sub-groups of $G=SU(2)$ or $SO(3)$ are finite, 1-dimensional,
   including $S^1$ and its normalizer $N(S^1)$, and $G$ itself. So $\phi^{-1}(0)$
   and $M_0$ may contain strata with some or all of these orbit types.

  We will remove strata of $M_0$ from $\phi^{-1}(U)/G$ in the order of lower dimensional ones
  (with ``bigger'' isotropy groups)
  to higher dimensional ones. By Lemma~\ref{lemstr}, we only need to
  prove that the link of each removed stratum is connected and
  simply connected.

    Assume there is an $(H)$-stratum in $M_0$.   By Theorem~\ref{thmform'}, a neighborhood $A$ in $M$ of a connected component
    of the
     $(H)$-stratum of $\phi^{-1}(0)$ is isomorphic to $A=G\times_H(\mathfrak{b}^{\bot}\times V)$, where $\mathfrak{b}^{\bot}$ is the
     annihilator of $\mathfrak{b}=Lie(H)$ in $\mathfrak{g}^*$.  Split $V=W\oplus V^H$.  By
   Lemma~\ref{lemlink}, the link of the connected component of the $(H)$-stratum in $A/G$ is $L_H=S(\mathfrak{b}^{\bot}\times W)/H$.

   1. Assume $H=G$. Then the above $L_H=L_G=S(W)/G$.
   Since the moment map image intersects a neighborhood of $0$,  $W$ has to be a
  nontrivial complex $G$-representation, so $S(W)$ (with high dimension) is connected and simply connected.   By Lemma~\ref{lembredon} or by Theorem~\ref{thmarm},
  $L_G$ is connected and simply connected.

   2. Assume $H=N(S^1)$. Then the above $L_H=S(\mathbb{R}^2\times W)/N(S^1)$, where $\mathbb{R}^2$ is such that
   $\mathbb{R}^2\oplus Lie^*(N(S^1))=\mathfrak g^*$. In this case, $N(S^1)$ acts on $\mathbb{R}^2$ as the action of
    $O(2)$. If $W=0$, $L_H$ is a point. Otherwise, by Theorem~\ref{thmarm}, $L_H$ is connected and simply connected.

   3. Assume $H=S^1$. Similar to 2.

   4. Assume $H=\Gamma$, where $\Gamma$ is a finite subgroup of $G$. Then the above link is
   $L_H=L_{\Gamma}=S(\mathbb{R}^3\times W)/\Gamma$.
   Since each element of $\Gamma$ belongs to a maximal torus
    of $G$,  it has a non-zero fixed vector in
    $S(\mathbb{R}^3)$. By Theorem~\ref{thmarm}, $L_{\Gamma}$ is
    connected and simply connected.
\end{proof}

   \section{Proof of Theorem 3}
    In this section, we prove Theorem 3.\\

    Let $G$ be a connected compact non-abelian Lie group with Lie algebra $\mathfrak g$ and dual Lie algebra $\mathfrak g^*$. Let
    $\mathfrak{t}_+^*$ be a fixed closed positive Weyl chamber of $\mathfrak g^*$. The set  $\mathfrak{t}_+^*$
    intersects each coadjoint orbit at a unique point, and it consists of open faces with different
    dimensions. All the points on a fixed open face of   $\mathfrak{t}_+^*$ have the same stabilizer
    group under the coadjoint action. We advise
    the reader to distinguish Weyl chamber and the faces of the Weyl chamber in this section with the $\mathbf{chamber}$
    and $\mathbf{faces}$ of the abelian moment polytope in Section 3 for which we used and we will use black letters.

    By Theorem~\ref{nonabelian-convexity},
    $\phi(M)\cap\mathfrak{t}_+^*=\triangle'$ is a convex polytope.
   Let us call the highest dimensional face $\tau^P$ of $\mathfrak{t}_+^*$ which contains values of $\phi$ the
   $\mathbf{principal\, face}$, and let us call the generic values of $\phi$ on $G\cdot\tau^P$ generic values of $\phi$.
   Let $G_{\tau^P}$ (connected) be the stabilizer group of $\tau^P$ under the coadjoint action.
   Let $U_{\tau^P}$ be the slice at $\tau^P$.  By Theorem~\ref{thmcross},
   $G_{\tau^P}$ acts on the $\mathbf{principal\, cross\, section}$
   $\phi^{-1}(U_{\tau^P})$.  Split $G_{\tau^P}=G'_{\tau^P}\times T_{\tau^P}$, where
    $G'_{\tau^P}$ is semi-simple, and $T_{\tau^P}$ is abelian. By Theorem 3.1 in \cite{LMTW} (see the cited theorem below),
    the semi-simple part $G'_{\tau^P}$
   acts trivially,
   only the connected central torus $T_{\tau^P}$ of $G_{\tau^P}$ acts on
   $\phi^{-1}(U_{\tau^P})$ non-trivially. If $\tau^P$
   is in the open positive Weyl chamber $\sigma$ of $\mathfrak{t}_+^*$, then
   $\phi^{-1}(U_{\sigma})$ has the maximal torus $T$ action.
   Otherwise, the
    central torus  $T_{\tau^P}$  of $G_{\tau^P}$ which acts on $\phi^{-1}(U_{\tau^P})$ has
    a smaller dimension than the dimension of the maximal torus $T$ of $G$.\\

    Let us write Theorem 3.1 in \cite{LMTW} in the following
    \begin{theorem}\label{thmLMTW}
    Let $G$ be a compact connected Lie group, and $M$ a connected Hamiltonian $G$-manifold with moment map
    $\phi: M\rightarrow\mathfrak{g}^*$.

    a. There exists a unique open face $\tau^P$ of the Weyl chamber $\mathfrak{t}^*_+$ with the property
       that $\phi(M)\cap\tau^P$ is dense in $\phi(M)\cap\mathfrak{t}^*_+$.

    b. The preimage $Y=\phi^{-1}(\tau^P)$ is a connected symplectic $T_{\tau^P}$-invariant submanifold of $M$, and
       the restriction $\phi|_Y$ of $\phi$ to $Y$ is a moment map for the action of $T_{\tau^P}$.

    c. The set $G\cdot Y=\{g\cdot m\,|\,g\in G, m\in Y\}$ is dense in $M$.
    \end{theorem}

    \begin{remark}\label{remstabilizer}
    For Theorem 3, without loss of generality, let us assume that the
    $\mathbf{principal\, stabilizer\, type}$ of the points in $M$
    intersects the center of $G$ trivially. If not, the above intersection
    subgroup (which is contained in the center) is contained in all
    the stabilizer groups of the points in $M$, so we can divide it
    out and consider the quotient group action.
    \end{remark}

    The proof of Theorem 3 consists of Lemma~\ref{lem6.1}, Lemma~\ref{lem6.2}, and  Lemma~\ref{lem6.3}.

    \begin{lemma}\label{lem6.1}
      Let $(M, \omega)$ be a connected, compact symplectic manifold equipped with a Hamiltonian  $G$ action
  with moment map $\phi$, where $G$
   is a connected compact non-abelian Lie group.
   Assume that a $G$-invariant metric is chosen on $\mathfrak{g}^*$.
   Let $b\in\mathfrak{t}_+^*$ be the furthest moment map value from the origin. Then $\pi_1(M)\cong\pi_1(M_b)$.
   \end{lemma}

   \begin{proof}
    Let the Weyl group act on $\triangle'$. The image of $\triangle'$ under
    this action is a polyhedron set (not necessarily convex) in $\mathfrak{t}^*$,
    and, we see that the
    Weyl group orbit through $b$ consists of the furthest points in $\mathfrak{t}^*$ to
    the origin, they are some boundary vertices of the polyhedron.

         The maximal torus $T$ of $G$ acts on $M$ with moment map $\phi_T$ being the
         projection to $\mathfrak{t}^*$ of the $G$ moment map $\phi$. By Theorem~\ref{abelian-convexity},
         this moment map image is
         a convex polytope.
         The point $b$ (and its Weyl group images)
         is a boundary vertex of the image of $\phi_T$. To see this, take
         the line segment $ob$, take the hyperplane in $\mathfrak{g}^*$ perpendicular to
         $ob$. Only the $G$ moment map images on this hyperplane
         will be projected to $b$. But, there cannot be points other
         than $b$ on this hyperplane which are in the $G$ moment map
         image. If there was, then this point is further than $b$ to
         the origin, and the intersection of $\mathfrak{t}^*_+$ and
         the coadjoint orbit through this point would be a point on
         $\triangle'$ which is further than $b$ to the origin. So we
         have proved that $\phi^{-1}(b)=\phi_T^{-1}(b)$ is a fixed
         point set component of the $T$ action on $M$. Therefore
         $\phi^{-1}(b)$ is a compact symplectic manifold with a
         trivial $T$ action.

       Let $G_b$ be the stabilizer group of $b$ under the coadjoint action. Then
       $G_b$ acts on $\phi^{-1}(b)$. The maximal torus of
       $G_b$ is also $T$. Since $T$ acts trivially on
       $\phi^{-1}(b)$, $G_b$ acts trivially on $\phi^{-1}(b)$.
       So $\pi_1(M_b)\cong\pi_1(\phi^{-1}(b))\cong\pi_1(M)$ by
       Lemma~\ref{lem3.1}.
     \end{proof}

   \begin{lemma}\label{lem6.2}
      Let $(M, \omega)$ be a connected, compact symplectic manifold equipped with a Hamiltonian  $G$ action
  with moment map $\phi$, where $G$
   is a connected compact non-abelian Lie group.
   Let $\tau^P\subset\mathfrak{t}_+^*$ be the $\mathbf{principal\,
   face}$, and let
    $a, b\in \tau^P$ be any two moment map values. Then $\pi_1(M_a)\cong\pi_1(M_b)$.
    \end{lemma}
    \begin{proof}
    By Theorem~\ref{thmLMTW}, the $\mathbf{principal\, cross\, section}$  $\phi^{-1}(U_{\tau^P})=Y$
    is a connected symplectic submanifold of $M$ with a torus $T_{\tau^P}$
    action whose moment map is $\phi|_Y$. Since $\phi|_Y$ is proper onto its
    image, we can use Theorem~\ref{thmLW}. So we still have
    Lemma~\ref{lem3.4} which implies Lemma~\ref{lem3.3}.
    By Lemma~\ref{lem3.2} and Lemma~\ref{lem3.3}, we have
    $\pi_1(Y_a)\cong\pi_1(Y_b)$ for all $a, b\in \tau^P$.
    By definition of the reduced spaces,
    this is to say $\pi_1(M_a)\cong\pi_1(M_b)$.
   \end{proof}

    Lemma~\ref{lem6.1} and Lemma~\ref{lem6.2} immediately imply the following special case of
     Theorem~\ref{thm3}.
    \begin{corollary}\label{cor6}
  Let $(M, \omega)$ be a connected, compact symplectic manifold equipped with a Hamiltonian  $G$ action
  with moment map $\phi$, where $G$
   is a connected compact non-abelian Lie group.
   Assume that the polytope $\phi(M)\cap\mathfrak{t}_+^*=\triangle'$ only lies on one face of the positive Weyl chamber,
   i.e., the moment polytope  only has a $\mathbf{principal\,
   face}$.
   Then, as fundamental groups
   of topological spaces, $\pi_1(M)\cong\pi_1(M_{red})$, where $M_{red}$ is the symplectic quotient at any coadjoint orbit of the moment map $\phi$.
   \end{corollary}

  For general cases, if we have Lemma~\ref{lem6.1} and Lemma~\ref{lem6.2}, we are only left  to
  show the following

    \begin{lemma}\label{lem6.3}
      Let $(M, \omega)$ be a connected compact symplectic manifold equipped with a Hamiltonian  $G$ action
  with moment map $\phi$, where $G$
   is a connected compact non-abelian Lie group. Let $\tau^P$ be the  $\mathbf{principal\,
   face}$ of the moment polytope.
   Let $c$ be a value of $\phi$ which is not on $\tau^P$, i.e.,
   $c$ is on a lower dimensional face of $\mathfrak{t}^*_+$.
     Let $a$ be a generic value on $\tau^P$ very near $c$. Then
   $\pi_1(M_c)\cong\pi_1(M_a)$.
   \end{lemma}

    Now, assume $c\in\tau$, where $\tau$ is a face of
    $\mathfrak{t}^*_+$. Let $G_{\tau}$ be the stabilizer group of
    points on $\tau$ under the coadjoint action. Then $G_{\tau}$ is a compact and connected Lie
    subgroup. Clearly, $G_{\tau}$ contains $T$. Let $U_c$ be the
    natural slice (see Section 4) at $c$. By Theorem~\ref{thmcross}, $U_c$ contains $\tau^P$. Let $R=\phi^{-1}(U_c)$ be
    the cross section (see Theorem~\ref{thmcross}) on which
    $G_{\tau}$ acts with moment map being $\phi|_R$. Since $G_{\tau}\cdot\phi^{-1}(\tau^P)$ is open, dense and connected
    in $R$ by Theorem~\ref{thmLMTW}, $R$ is connected.
     By definition of
     the reduced spaces at $c$ and at $a$, to compare $\pi_1$ of the
     two quotients $M_a$ and $M_c$, we may restrict attention to $(R,
     \omega|_R,
     G_{\tau}, \phi|_R)$.
     Split
    $G_{\tau}=G_1\times T_2$, where $G_1$ is semi-simple, and
    $T_2$ is abelian and it is the connected component of the center of $G_{\tau}$.
    The linear space spanned by the
     face $\tau$ is the dual Lie algebra of $T_2$.
     Since $\tau$ lies on the center of $\mathfrak{g}^*_{\tau}$, we may
     shift the moment map $\phi|_R$ by $c$ such that the value $c$ corresponds to the value $0$ of the shifted
      moment map  $\phi'|_R$. So, without loss of generality, we assume that we
      have the Hamiltonian space $(R, \omega|_R,
     G_{\tau}, \phi'|_R)$, and we want to prove that
     $\pi_1(R_0)\cong\pi_1(R_{G_{\tau}\cdot a})$, where $a\in\tau^P$ is a generic value  near $0$ in
     $\mathfrak{g}^*_{\tau}$. Although $R$ may not be compact, since $M$ is compact,
     the moment map $\phi|_R$ is proper onto its image. So, if we
     take a small $G_{\tau}$-invariant neighborhood $U$ of $0$ in $\mathfrak{g}^*_{\tau}$,
     we may assume that the $G$-equivariant gradient flow
     of $f=\|\phi\|^2$ restricts to the $G_{\tau}$-equivariant gradient flow of $\|\phi|_R\|^2$ on
       $\phi|_R^{-1}(U)$.
     By Theorem~\ref{thmLW}, there exists a smaller $G_{\tau}$-invariant neighborhood $U'$ of $0$, such that
     $\phi|_R^{-1}(U')$ is $G_{\tau}$-equivariantly homotopy equivalent to $\phi|_R^{-1}(0)$. \\

      So we only need to prove the following general lemma

   \begin{lemma}\label{lem6.4}
      Let $(N, \omega)$ be a connected  symplectic manifold equipped with a Hamiltonian  $K$ action
  with moment map $\phi$ which is proper onto its image, where $K$
   is a connected compact non-abelian Lie group.
   Assume that the moment map takes value at $0$ and it takes values in a neighborhood of
   $0$ in $\mathfrak{k}^*$ (not necessarily that the moment map image fills an open neighborhood
    of $0$).
   Let $a$ be a generic value near $0$ in the positive Weyl chamber $\mathfrak{t}^*_+$.
   Then
   $\pi_1(N_0)\cong\pi_1(N_a)$.
   \end{lemma}

   \begin{lemma}\label{lem6.5}
      Under the assumptions of Lemma~\ref{lem6.4}, assume in
      addition that
     $U\subset\mathfrak{k}^*$ is a small open invariant neighborhood of $0$
    such that $\phi^{-1}(U)$ equivariantly deformation retracts
    to $\phi^{-1}(0)$.
    Let $B$ be the set of values in $U$ which are on the faces other than the $\mathbf{principal\, face}$ $\tau^P$
    of the closed
    positive Weyl chamber and the set of  values which are on $\tau^P$ but not on the open connected
    $\mathbf{chamber}$ of generic values containing $a$.
   Then
   $\pi_1(\phi^{-1}(U)/K)\cong\pi_1(\phi^{-1}(U)/K-\phi^{-1}(K\cdot B)/K)$.
   \end{lemma}

   Assume we have this lemma. Then the proof of Lemma~\ref{lem6.4}
   goes as the following.
   \begin{proof}
   Since  $\phi^{-1}(U)$ equivariantly deformation retracts
      to $\phi^{-1}(0)$,
      we have $\pi_1(\phi^{-1}(U)/K)\cong\pi_1(N_0)$.
   By Lemma~\ref{lem6.5},
   $\pi_1(\phi^{-1}(U)/K)\cong\pi_1(\phi^{-1}(U)/K-\phi^{-1}(K\cdot B)/K)\cong\pi_1(N_a)$.
   \end{proof}

   It is left to prove Lemma~\ref{lem6.5}. The following two facts about the coadjoint action of a connected compact
   Lie group are needed in the proof of Lemma~\ref{lem6.5}.
    \begin{proposition}\label{propfix}
       Let $G$ be an n-dimensional connected compact Lie group with Lie algebra
    $\mathfrak{g}$ and dual Lie algebra $\mathfrak{g}^*$.
    Let $H\neq G$ be a subgroup of $G$  with Lie algebra $Lie(H)=\mathfrak b$,
    and let $\mathfrak b^{\circ}$ be the annihilator of $\mathfrak
    b$ in $\mathfrak{g}^*$. The space $\mathfrak b^{\circ}$ can be identified with the
    orthogonal  complement $\mathfrak b^{\bot}$ of $\mathfrak b$ in $\mathfrak{g}$ for a suitable metric.
    The subgroup $H$ acts on $\mathfrak b^{\circ}$ by
    the coadjoint action.
   Then, the smallest normal subgroup $N_H$ of $H$ containing the
    identity component of
     $H$ and all those elements of $H$ which have non-zero fixed points is $H$ itself.
     \end{proposition}

    \begin{proof}
   (1) If $H$ is connected, we are done.

   (2) Assume that $H$ is not connected.
       Let $H^0$ be the identity component of $H$. Then $H$ is generated by $H^0$ and
       finitely many elements, say $h_1, h_2, ..., h_k$. If we can prove that each
       $h_i, i=1,..., k$ has a non-zero fixed point in $\mathfrak b^{\circ}$, then we are done.
     Let $T_1$ be a maximal torus of $H$. Then $T_1\subset H^0$.
     Let $T$ be a  maximal torus of $G$ such that $T=T_1\times T_2$, where the dual Lie algebra $\mathfrak t_2$ of
     $T_2$ is contained in  $\mathfrak b^{\circ}$ ($T_1$ or $T_2$ can be trivial).
     If $h_i\in T$, then $\mathfrak t_2\neq 0$. So $h_i$ fixes the subspace $\mathfrak t_2$ of $\mathfrak b^{\circ}$.
     If $h_i\notin T$, i.e., $h_i$ is in a different maximal torus $T'$ other than $T$, then the dual
         Lie algebra $\mathfrak{t}'^*$ of $T'$ has to have a non-zero component in $\mathfrak b^{\circ}$. Indeed,
     if $\mathfrak{t}'^*\subset\mathfrak b$, then $h_i\in T'\subset H^0$, a contradiction. So $h_i$ fixes the above
     non-zero component in $\mathfrak b^{\circ}$.
      \end{proof}

 \begin{proposition}\label{propLie}
    Let $G$ be a connected compact semi-simple non-abelian Lie group with Lie algebra $\mathfrak{g}$
    and dual Lie algebra $\mathfrak{g}^*$.
  Let $H$ be a subgroup of $G$ with
    Lie algebra $\mathfrak{h}$. Let $\mathfrak{h}^{\circ}$ be the annihilator of $\mathfrak{h}$ in
$\mathfrak{g}^*$. If
     $H\neq G$, then $H$ has at least codimension 2. If $H$ has codimension 2, i.e.,
     if dim $(\mathfrak{h}^{\circ})=2$, then $H$ acts
     (coadjoint action)
    on $S(\mathfrak{h}^{\circ})=S^1$ transitively.
   \end{proposition}
   \begin{proof}
    Let $\mathfrak g=\mathfrak t\oplus \bigoplus_{\alpha\in\mathcal
    R^+}M_{\alpha}$ be the (real) root space decomposition of $\mathfrak
    g$, where $\mathfrak t$ is the Lie algebra of the maximal torus $T_G$ of $G$, and
    $\mathcal R^+$ is the set of positive roots. We know that each
    root space $M_{\alpha}$ is 2-dimensional, and for two generators $X, Y\in M_{\alpha}$, there is a $Z\in \mathfrak t$, such that
    $[Z, X]=2X, [Z, Y]=-2Y, [X, Y]=2Z$ (this corresponds to the Lie algebra of $SU(2)$). Let $T_H$ be
    the maximal torus  of $H$. We can split $\mathfrak{h}$ similar to the splitting of $\mathfrak g$.
     If dim($T_H)<$dim$(T_G)$, then at least one vector of $\mathfrak t$
    is missing in $Lie(T_H)\subset\mathfrak{h}$, and at least one $M_{\alpha}$ is
    missing in $\mathfrak{h}$. Indeed, let us assume that dim($T_H)=$dim$(T_G)-1$. Let $Z$ be the non-zero
    vector in $\mathfrak t$ but not in $Lie (T_H)$. It is clear that the linear space $M_{\alpha}$ which has the above mentioned
    property with $Z$ cannot be in $\mathfrak{h}$ since
    $\mathfrak{h}$ is a Lie algebra.
    If only one generator $X\in M_{\alpha}$ were in $\mathfrak{h}$, then $X$ would contribute to
    a vector in $Lie(T_H)$ which contradicts to the fact that dim($T_H)=$dim$(T_G)-1$. Therefore,
     in the case that    dim($T_H)<$dim$(T_G)$,
    $H$ has at least codimension 3 in
    $G$. If  dim($T_H)=$dim$(T_G)$, then $\mathfrak t\subset\mathfrak{h}$. Since $H\neq G$, at least one
    $M_{\alpha}$ is missing in $\mathfrak{h}$. So, in this case, $H$
    has at least codimension 2 in $G$. When codimension of $H$ is 2,
    $\mathfrak{h}^{\circ}=M_{\alpha}$ for some $\alpha$. So
    there is an $S^1$ (generated by the above $Z$) in the maximal torus of $G$ which acts on
    $S(\mathfrak{h}^{\circ})=S^1$ transitively.
    \end{proof}

   Let us now proceed to prove Lemma~\ref{lem6.5}. For the same reason as we made in Remark~\ref{order},
   we will prove Lemma~\ref{lem6.5} by  induction on removing
   $\phi^{-1}(K\cdot\tau)/K$ in the order of lower dimensional faces ($\tau$s' whose preimage may have
   bigger stabilizer groups) to higher dimensional
   faces of the closed positive Weyl chamber. Let $C$ be the central face of the closed
   positive Weyl chamber
   $\mathfrak{t}^*_+$. Write
$K=K_1\times T_1$,
   where $K_1$ is semi-simple, and $T_1$ is abelian (both $K_1$ and $T_1$ are connected). If $C$ is the
   $\mathbf{principal\, face}$, then by Theorem~\ref{thmLMTW}, only the central torus $T_1$ acts on $N$.
   Corollary~\ref{cor6} addressed this case.
   In the following, we assume
    that $C$ is not the only face which contains the image of $\phi$. Let us first remove $\phi^{-1}(C)/K$ from
    $\phi^{-1}(U)/K$.

     \begin{lemma}\label{lem6.6}
     Under the assumptions of Lemma~\ref{lem6.5}, let $C$ be the central face of $\mathfrak{t}^*_+$.
     Assume that
       $C$ is not the only face which contains the image of $\phi$. Then
      $\pi_1(\phi^{-1}(U)/K)\cong\pi_1(\phi^{-1}(U)/K-\phi^{-1}(C)/K)$.
     \end{lemma}
     \begin{proof}
        Let $K=K_1\times T_1$, where $K_1$ is semi-simple, and
         $T_1$ is abelian. Both $K_1$ and $T_1$ are connected by assumption.
     Let the Lie algebra of $K_1$ be $\mathfrak{k}_1$, and let its dual Lie algebra be
         $\mathfrak{k}_1^*$. Let the Lie algebra of $T_1$ be $\mathfrak{t}_1$,
         and let its dual Lie algebra  be $\mathfrak{t}^*_1$.  Then $C=\mathfrak{t}^*_1$.

       We will remove
      strata of  $\phi^{-1}(C)/K$ in the order of lower dimensional
    ones to higher dimensional ones and use Lemma~\ref{lemstr} repeatedly.

     Now, assume there is a stratum with isotropy type $(H)$ (subgroups conjugate to $H$) in  $\phi^{-1}(C)/K$.
     Then $H$ has the form of a product $H=H_1\times T'$, where $H_1$ is a
     subgroup of $K_1$, and $T'$ is a subgroup of $T_1$. Let $\mathfrak{h}_1^{\bot}$ be the annihilator of the Lie
     algebra $\mathfrak{h}_1=Lie(H_1)$ in $\mathfrak{k}_1^*$, and let $(\mathfrak{t}')^{\bot}$
      be the annihilator of the Lie algebra
     $\mathfrak{t}'=Lie(T')$ in $\mathfrak{t}^*_1$.  Since every orbit in  $\phi^{-1}(C)$
     is isotropic, we can use Theorem~\ref{thmform'}.
     By Theorem~\ref{thmform'}, a
     neighborhood in $N$ of an orbit in  $\phi^{-1}(C)$ with isotropy type
     $(H)$ is isomorphic to $A=K\times_H (\mathfrak{b}^{\bot}\times V)$, where $\mathfrak{b}^{\bot}$ is the annihilator of
     $\mathfrak{b}=Lie(H)$ in $\mathfrak{k}^*=Lie^*(K)$.
     Let $V=W\times V^H$. Write $A=K\times_H (\mathfrak{h}_1^{\bot}\times(\mathfrak{t}')^{\bot}\times (W\times V^H))$.
     The quotient $A/K$ is $(\mathfrak{h}_1^{\bot}\times W)/H\times(\mathfrak{t}')^{\bot}\times V^H$.
     The $(H)$-stratum of $\phi^{-1}(C)/K$ in $A/K$ is
      $(\mathfrak{t}')^{\bot}\times V^H$.
      The link $L_H$ of this connected $(H)$-stratum is
      $S(\mathfrak{h}_1^{\bot}\times W)/H$. Now, we consider all the possible cases of $H$.

     1.  Assume that $H_1=K_1$. Then the above link is
    $S(W)/(K_1\times T')$.
           Since we assumed that
    the moment map image intersects not only the central face of the positive Weyl chamber,
    $W$ has to be a non-trivial complex $K_1$ representation. So $S(W)$ (with dimension at least 3) is connected and simply connected. Due to the fact that
  the quotient of $W$ by a finite abelian group $\Gamma$ is homeomorphic
 to $W$, we may assume that $T'$ is connected.
   By Lemma~\ref{lembredon}, the link
      $L_H$
     is connected and simply connected.

    2. Assume $H_1\neq K_1$.  By Proposition~\ref{propLie}, dim$(\mathfrak{h}_1^{\bot})\geq 2$.

      First, let us assume that $W\neq 0$, then dim$(\mathfrak{h}_1^{\bot}\times W)>2$. So $S(\mathfrak{h}_1^{\bot}\times
    W)$ is connected and simply connected.
    By Proposition~\ref{propfix} and Theorem~\ref{thmarm} ($T'$ acts on $\mathfrak{h}_1^{\bot}$ trivially), the link
    $S(\mathfrak{h}_1^{\bot}\times W)/(H_1\times T')$
    is connected and simply connected.

    Next, let us assume that dim$(\mathfrak{h}_1^{\bot})=2$ and $W=0$. The link is $S(\mathfrak{h}_1^{\bot})/H$.
      By Proposition~\ref{propLie}, this link is a point, therefore connected and simply connected.

    Now, let us  assume that $\mathfrak{h}_1^{\bot}>2$ and $W=0$.
    By Proposition~\ref{propfix} and Theorem~\ref{thmarm}, the link is connected and simply
    connected.
    \end{proof}
    Proof of Lemma~\ref{lem6.5}:
   \begin{proof}
    By assumption, $0\in im(\phi)$ and $im(\phi)$ intersects a small neighborhood of $0$. If the image of $\phi$ only lies
    on the central face $C$ of the positive Weyl chamber, then the main theorem comes down
    to the case of Corollary~\ref{cor6}.
    Now, we assume that the image of $\phi$ intersects not only one face of the positive Weyl chamber.

     Lemma~\ref{lem6.6}  removed  $\phi^{-1}(C)/K$ from $\phi^{-1}(U)/K$.

    Assume now that the image of $\phi$ intersects another  higher dimensional face $\tau$ other than $\tau^P$ of the positive
    Weyl chamber. Suppose the stabilizer group of $\tau$ under the
    coadjoint action is $K_{\tau}$. Let $U_{\tau}$ be the natural slice at $\tau$. Then by the symplectic cross
    section theorem, $Y_{\tau}=\phi^{-1}(U_{\tau})$  is a symplectic
    submanifold with a $K_{\tau}$ action. The face $\tau$
    lies on the central dual Lie algebra of $K_{\tau}$.
    Similar to  Lemma~\ref{lem6.6}, we remove
    $\phi^{-1}(\tau)/K_{\tau}$ from $(\phi^{-1}(U)\cap
    Y_{\tau})/K_{\tau}$. In order to  remove
    $\phi^{-1}(K\cdot\tau)/K$ from $\phi^{-1}(U)/K$, we only need to notice that, by
    equivariance of the moment map, each stratum of
    $\phi^{-1}(K\cdot\tau)/K$ has the same link in $\phi^{-1}(U)/K$
    as the corresponding stratum of  $\phi^{-1}(\tau)/K_{\tau}$ in $(\phi^{-1}(U)\cap
    Y_{\tau})/K_{\tau}$. Indeed, to prove this, observe that $K\cdot Y_{\tau}$ is equivariantly
    diffeomorphic to the bundle $K\times_{K_{\tau}} Y_{\tau}$ over the coadjoint orbit $K/K_{\tau}$.
    So each connected stratum $\tilde S$ in  $K\cdot Y_{\tau}$ is $\tilde S=K\times_{K_{\tau}}\tilde S'$, where
    $\tilde S'$ is the corresponding stratum in $Y_{\tau}$. A neighborhood of $\tilde S$ in $K\cdot Y_{\tau}$
    corresponds to
    a neighborhhood of $\tilde S'$ in $Y_{\tau}$.  Therefore the quotient stratum $\tilde S/K$ in $K\cdot Y_{\tau}/K$
    has the same link as $\tilde S'/K_{\tau}$ in $Y_{\tau}/K_{\tau}$. If we restrict to the invariant set
    $\phi^{-1}(U)$, the same property still holds.

    We can remove similarly  $\phi^{-1}(K\cdot\tau')/K$ from $\phi^{-1}(U)/K$
    if the image of $\phi$ intersects other non-principal faces $\tau'$s' of the positive Weyl chamber.

     Now assume only   $\phi^{-1}(K\cdot(\tau^P\cap U))/K$ is
     remaining. If all the values on $\tau^P\cap U$ are regular, we
     are done. Otherwise, we will remove and ``flow" in the
     following way untill only the connected open $\mathbf{chamber}$
     containing $a$
     is remaining. After we ``removed" all the non-principal faces
     of the positive Weyl chamber, if on the ``verge" of $\tau^P\cap
     U$, there is an open connected $\mathbf{chamber}$ $U'$ not containing $a$, we use the
     gradient flow of suitable components of the $T$ moment map to deformation retract
     $\phi^{-1}(U')$ to  $\phi^{-1}(\mathcal{F}s)$, where $\mathcal{F}s$ are
     certain singular $\mathbf{faces}$ around $U'$. Correspondingly, by equivariance of $\phi$ again,
     $\phi^{-1}(K\cdot U')$ deformation retracts to  $\phi^{-1}(K\cdot \mathcal{F}s)$.
     Then, we use
     Lemma~\ref{lem3.4} to remove $\phi^{-1}(
     \mathcal{F}s)/T_{\tau^P}$ from the remaining part of
      $\phi^{-1}(\tau^P\cap U)/T_{\tau^P}$,
      and we
     use equivariance as above to remove $\phi^{-1}(K\cdot \mathcal{F}s)/K$.
     Or, we may only need to do removing if there is only one connected open $\mathbf{chamber}$ which is the one containing $a$ is left
     in $\tau^P\cap U$ (the other singular $\mathbf{faces}$ are in the closure of this $\mathbf{chamber}$).
     We may need to repeat the procedure  untill
    only the connected open $\mathbf{chamber}$ containing $a$
    remains.

     Now, we have removed all $\phi^{-1}(K\cdot B)/K$ from $\phi^{-1}(U)/K$. The
     remaining space
    is homotopy equivalent to $\phi^{-1}(K\cdot a)/K=N_{K\cdot a}$, where $a$ is a generic value
    in the $\mathbf{chamber}$ containing $a$.
    \end{proof}
   \section{proof of Theorem~\ref{thm4}}

    In this section, we will prove Theorem~\ref{thm4}. For simply connected manifold $M$, Armstrong's theorem
    (Theorem~\ref{thmarm}) tells us that $\pi_1(M)\cong\pi_1(M/G)$.
     In this case, by Theorem~\ref{thm1}, Theorem~\ref{thm2} and Theorem~\ref{thm3}, we have Theorem~\ref{thm4}. In this section,
     we use our method of removing to give a direct proof of Theorem~\ref{thm4}.
    Since we proved that all reduced spaces have isomorphic fundamental groups, we only need to prove
    that $M/G$ has the same fundamental group as that of a particular reduced space. \\

    In Lemma~\ref{lem3.4}, we did removing from ``one side'', i.e., we took $\bar{U'}$. Now, let us do a removing from
    ``all sides'', i.e., we take $U$ itself. The proof is even simpler.
        \begin{lemma}\label{lem3.4'}
      Assume we have the assumptions of Theorem~\ref{thm4}, where $G$ is a torus $T$.
      Let $\mathcal{F}$ be a $\mathbf{face}$ of the moment polytope $\triangle$ which is not a $\mathbf{chamber}$.
      Let $U$ be a small open neighborhood of $\mathcal{F}$ on $\triangle$ ($U$ does not intersect the
      $\mathbf{faces}$ which are in the closure of $\mathcal{F}$).
      Let $S$ be the set of singular orbits (or non-generic orbits) in $\phi^{-1}(\mathcal{F})$.
      Then
      $\pi_1(\phi^{-1}(U)/T)\cong\pi_1(\phi^{-1}(U)/T-S/T)$.
      \end{lemma}

         \begin{proof}
     Assume that the dimension of the face $\mathcal F$ is $m$ with $m\geq 0$. We do the removing
     from lower dimensional strata to higher dimensional strata. As we did before, we only need to
     check that the link of the removed stratum is connected and simply connected.

     As in the proof of Lemma~\ref{lem3.4}, we split the torus $T=T^{n-m}\times T^m$.
     The possible
     stabilizer types of the strata in $\phi^{-1}(\mathcal F)$
     have the form $H=(T_1\times\Gamma')\times\Gamma$, where $T_1$ is a connected subgroup
     of  $T^{n-m}$, $\Gamma'$ is a finite subgroup of $T^{n-m}$, and $\Gamma$ is a finite subgroup of $T^m$.

        By Theorem~\ref{thmform}, a
        neighborhood in $M$ of an orbit with
        isotropy type $(H)$ is isomorphic to
        $A=T\times_H(\mathbb{R}^l\times\mathbb{R}^m\times V)$. Split $V=W\times V^H$.
         If $W=0$, then all
    the points in $A$ have the same stabilizer group. Since $A$ is open, this stabilizer group is the generic stabilizer group.
    So we assume that $W\neq 0$ by our assumption.
        The $(H)$-stratum in $A$ which was mapped to $\mathcal F$ is
        $T\times_H(\mathbb{R}^m\times V^H)$.
    The quotient of $A$ by $T$ is
    $A/T=\mathbb{R}^l\times\mathbb{R}^m\times V^H\times W/H$.
    The link $L_H$ of the corresponding quotient $(H)$-stratum
        in $A/T$ is $S(\mathbb{R}^l\times W)/H$. The possible cases of $H$ (or of $l$) are:

    1. In the case $l=0$,
       since the moment map value fills the neighborhood $U$, $W$ has to be a non-trivial
    representation of $H$ with a non-trivial moment map. So either $L_H$ is a point, in the case of dim $(W)=2$,
     therefore connected and simply connected,
    or, it is connected and
        simply connected by the fact that $W/\Gamma$ is homeomorphic to $W$ and by applying
        Lemma~\ref{lembredon} .

       2. In the case $l\neq 0$,
         $S(\mathbb{R}^l\times W)$ is simply connected ($W$ is a non-trivial
    complex vector space). By
        Theorem~\ref{thmarm}, $L_H$ is connected and simply
               connected.
     \end{proof}

   Lemma~\ref{lem3.4}, Lemma~\ref{lem3.4'} and  Lemma~\ref{lem6.6} did ``local removing'', namely, we chose
      small neighborhoods of a value or of a $\mathbf{face}$ on the moment polytope. Note
      that the removing itself does not depend on the gradient flow.
      In  Lemma~\ref{lem3.4} and Lemma~\ref{lem6.6}, we chose $U$
    to be small neighborhoods of one value. As in
    Lemma~\ref{lem3.4'}, if only for the purpose of removing (not
    for the purpose of comparing $\pi_1$ of the reduced spaces at
    nearby values), we could have taken $U$ in Lemma~\ref{lem3.4} and Lemma~\ref{lem6.6}
    to be an open small
    neighborhood of the $\mathbf{face}$ we are considering.
    Now, if we consider the global quotient $M/G$, we can do the
    same  removing in $M/G$, since the quotient of an orbit is only
    ``linked'' to its neighborhood in the space $M/G$.

      Another observation is that the gradient flow (or negative gradient flow) of the components of the moment map
      or of the moment map square  always retracts less singular regions on the manifold to more singular ones.

     Lemma~\ref{lem3.4}, Lemma~\ref{lem3.4'} and  Lemma~\ref{lem6.6}
     allow us to remove (the quotients of certain orbits). By using the gradient
     flow of suitable components of the moment map, we can retract
     regular regions to singular ones. These two operations are
     the main points of the proof of Theorem~\ref{thm4}. \\

      Proof of Theorem~\ref{thm4} for $S^1$ actions:

      \begin{proof}
    Assume that
   the moment map takes critical values at $a_0, a_1, ..., a_n$, and we have $a_0<a_1<...<a_n$.

  By Lemma~\ref{lem3.4} and an argument by the Van-Kampen theorem, we have
   $\pi_1(M/S^1)\cong\pi_1(M/S^1-M_{a_n})$. The right hand side is equal to $\pi_1(\phi^{-1}([a_0, a_n))/S^1)$.
   Now the negative gradient flow of $\phi$ equivariantly deformation retracts
   $\phi^{-1}([a_0, a_n))$ to $\phi^{-1}([a_0, a_{n-1}])$. So
   $\pi_1(\phi^{-1}([a_0, a_n))/S^1)\cong\pi_1(\phi^{-1}([a_0, a_{n-1}])/S^1)$.
   By Lemma~\ref{lem3.4} and a similar argument by the Van-Kampen theorem as the above,
   $\pi_1(\phi^{-1}([a_0, a_{n-1}])/S^1)\cong\pi_1(\phi^{-1}([a_0, a_{n-1}))/S^1)$. We use the
   negative gradient flow of $\phi$ again to ``flow down'' the next regular region $\phi^{-1}((a_{n-2}, a_{n-1}))$ to
   $\phi^{-1}([a_0, a_{n-2}])$.
   We repeatedly use the above procedure of
   removing and ``flowing'' until we reach $\pi_1(M/S^1)\cong\pi_1(M_{a_0})$. \\
    \end{proof}

    Proof of Theorem~\ref{thm4} for $G=T$ actions:

     \begin{proof}
     Different removing and flowing process can achieve the proof.
     The moment polytope $\triangle$ consists of $\mathbf{faces}$ with different
    dimensions; and, it may have one or more than one connected open
$\mathbf{chambers}$. The main
    tools we can use, or the main points of the proof are:
    (a) We can use  Lemma~\ref{lem3.4'} to remove singular
    $\mathbf{faces}$ on the boundary of $\triangle$.
    (b) The inverse image of a singular $\mathbf{face}$
     in the ``interior'' of $\triangle$ may contain
    regular orbits (orbits with generic stabilizer group). Whenever after we use Lemma~\ref{lem3.4'} for
    such a singular
    $\mathbf{face}$, the gradient flow of some components of the moment map
     always takes the remaining regular orbits to an appropriate regular region (the inverse image of an open
    $\mathbf{chamber}$).
    (c) We can start to remove from a vertex on the boundary of
    $\triangle$. We follow the principle we made in Remark~\ref{order} about the order of
    removing. We deform when it is allowed by using the gradient flow of
    suitable
    components of the moment map. When we encounter a removing of a
    singular $\mathbf{face}$ from the closure of only one connected
    $\mathbf{chamber}$, we may use Lemma~\ref{lem3.4}.
    (d) We can choose different process. In the end, we arrive at
    $\pi_1(M/G)\cong\pi_1(M_b)$, where $b$ is a value on $\triangle$.
    \end{proof}

  Proof of Theorem~\ref{thm4} for non-abelian $G$ actions:

      \begin{proof}
     We use the method of Lemma~\ref{lem6.6} to inductively remove $\phi^{-1}(G\cdot\tau)/G$ from
     $M/G$ for the
     faces $\tau$s' which are not the $\mathbf{principal\, face}$ $\tau^P$ of the closed positive Weyl chamber
     (Again, we may need the cross section theorem to do some removing in
     the cross section, and then use equivariance of the moment
     map. For this, see the proof of Lemma~\ref{lem6.5}. I would like to stress that,
     the removing itself does not depend on the gradient flow.).
     Assume, now, we have $\pi_1(M/G)\cong\pi_1(\phi^{-1}(G\cdot\tau^P)/G)$. By equivariance of the moment map,
     the two spaces $\phi^{-1}(G\cdot\tau^P)/G$ and $\phi^{-1}(\tau^P)/T$ are the same.
     We use the method for torus actions to prove that
     $\pi_1(\phi^{-1}(\tau^P)/T)\cong\pi_1(M_b)$ for some value $b\in\tau^P$.
  \end{proof}

\end{document}